\documentclass[12pt]{article}
\title{
Center stable manifolds around line solitary waves of the Zakharov--Kuznetsov equation with critical speed
\amssubj{35B35, 35Q53.}
\keywords{center stable manifolds; Zakharov--Kuznetsov equation; line solitary wave; transverse instability}
}
\author{YOHEI YAMAZAKI\footnote{{\it E-mail addresses:} yohei-yamazaki@hiroshima-u.ac.jp}\\  {\footnotesize Hirhoshima University,
 } \\  {\footnotesize  1-3-2 Kagamiyama, Higashi-Hiroshima City, Hiroshima, 739-8511 Japan} }
\usepackage{amssymb}
\usepackage{mathrsfs}
\usepackage{amsthm}
\usepackage{amsmath}
\usepackage{bm}
\usepackage{titlefoot}
\usepackage{fancyhdr}
\def\pdfliteral #1 {}

\numberwithin{equation}{section}

\newtheorem{theorem}{Theorem}[section]
\newtheorem{corollary}[theorem]{Corollary}

\newtheorem{lemma}[theorem]{Lemma}
\newtheorem{proposition}[theorem]{Proposition}

\theoremstyle{definition}
\newtheorem{definition}[theorem]{Definition}
\newtheorem{remark}[theorem]{Remark}

\renewcommand{\eqref}[1]{(\ref{#1})}

\renewcommand{\bigskip}{\vspace{0.3cm}}

\newcommand{\R}{{\mathbb R}}

\newcommand{\Z}{{\mathbb Z}}
\newcommand{\T}{{\mathbb T}}

\newcommand{\norm}[1]{{\left \lVert #1 \right \rVert}}

\newcommand{\tbr}[1]{\langle #1 \rangle}
\newcommand{\gbr}[1]{\lceil #1 \rfloor}

\newcommand{\RTL}{\mathbb {R} \times \mathbb {T}_L}

\date{}
\begin{document}

\maketitle

%{\small Keywords and phrases: center stable manifolds; Zakharov--Kuznetsov equation; line solitary wave; transverse instability}

\begin{abstract}
In this paper, we construct  center stable manifolds around unstable line solitary waves of the Zakharov--Kuznetsov equation on two dimensional cylindrical spaces $\RTL$ ($\T_L=\R/2\pi L\Z$).
In the paper \cite{YY5}, center stable manifolds around unstable line solitary waves have been constructed without critical speed $c \in \{4n^2/5L^2;n \in \Z, n>1\}$.
Since the linearized operator around line solitary waves with critical speed is degenerate, we prove the stability condition of the center stable manifold for critical speed by applying to the estimate of 4th order term of a Lyapunov function in \cite{YY2} and \cite{YY4}.
\end{abstract}

%{\small 2010 Mathematics Subject Classification.  35B35, 35Q53}

\section{Introduction}

We consider the two dimensional Zakharov--Kuznetsov equation
\begin{equation}\label{ZKeq}
 u_t + \partial_x(\Delta u + u^2)=0, \quad (t,x,y) \in \R\times \RTL,
\end{equation}
where $\Delta=\partial_x^2 + \partial_y^2$, $u=u(t,x,y)$ is an unknown real-valued function, $\T_L=\R/2\pi L \Z$ and $L>0$.
The equation \eqref{ZKeq} preserves the mass and the energy:
\[%\label{moment}
M(u) = \int_{\RTL} |u|^2 dxdy,
\]
\[%\label{energy}
E(u)= \int_{\RTL}\Bigl( \frac{1}{2}|\nabla u|^2 - \frac{1}{3}u^3 \Bigr) dxdy,
\]
where $u \in H^1(\RTL)$.

The Zakharov--Kuznetsov equation was introduced in \cite{Z K} to describe the propagation of inonic-acoustic waves in uniformly magnetized plasma.
The rigorous derivation of the Zakharov--Kuznetsov equation was proved in \cite{L L S}.
The Cauchy problem for the Zakharov-Kuznetsov equation is extensively studied in the literature.
The global well-posedness of the Zakharov--Kuznetsov equation in $H^s(\RTL)$ for $s>\frac{3}{2}$ has been proved by Linares, Pastor and Saut \cite{L P S} to study of the transverse instability of the N-soliton of the Korteweg-de Vries equation.
Molinet and Pilod \cite{M P} showed the global well-posedness in $H^1(\RTL)$ by proving a bilinear estimate in the context of Bourgain's spaces $X^{s,b}$.
For more results of the Cauchy problem of the Zakharov--Kuznetsov equation on whole spaces or torus, we refer to the papers \cite{AVF,AG1,G H,Ki,L P R T,L P 1,L P 2,L S,M S T}, and the references therein. 

The Zakharov--Kuznetsov equation has solitary wave solutions $\varphi(x-ct,y)$ which is nontrivial solutions to the stationary equation
\begin{equation}\label{eq-st}
-\Delta \varphi+c\varphi-\varphi^2=0
\end{equation}
for $c>0$.
de Bouard \cite{AdB} proved the orbital stability of positive solitary waves in $H^1(\R^2)$ which are ground state defined by the action 
\[S_c(u)=E(u)+\frac{c}{2}M(u)\]
of \eqref{eq-st} on $\R^2$.
C\^ote, Mu\~noz, Pilod and Simpson \cite{C M P S} proved the asymptotic stability of the positive solitary waves and multi solitary wave in $H^1(\R^2)$ by adapting the argument of Martel and Merle \cite{M M 1,M M 2,M M 3} to a multidimensional model.

The Zakharov--Kuznetsov equation has the line solitary wave
\[
Q_c(x-ct)=\frac{3c}{2}\cosh^{-2}\Bigl( \frac{\sqrt{c}(x-ct)}{2}\Bigr), \quad c>0
\]
which is also the one soliton of the Korteweg--de Vries equation.
It is well-known that the orbital stability of the one soliton $Q_c(x-ct)$ to the KdV equation was proved by  Benjamin \cite{TBB}.
The asymptotic stability of the one soliton on the exponentially weighted space was shown by Pego and Weinstein \cite{P W 2}.
This result was refined by Mizumachi \cite{TM1}, who treated perturbations in polynomially weighted spaces.
In \cite{M M 1,M M 2,M M 3}, Martel and Merle proved the asymptotic stability of the one soliton to the KdV equation on the energy space by using the Liouville type theorem and the monotonicity property.

The instability of the line solitary waves of the Zakharov--Kuznetsov equation on $\R^2$ was proved by Rousset and Tzvetkov \cite{R T 0}.
Johnson \cite{MAJ} proved the linear instability of line periodic solitary waves of  the Zakharov--Kuznetsov equation on $\T_{L_1} \times \T_{L_2}$ with sufficiently large $L_2$ by applying Evans function method.
Bridges \cite{TJB} showed the instability of the line solitary waves $Q_c(x-ct)$ of the Zakharov--Kuznetsov equation on $\RTL$ with sufficiently large traveling speed $c$.
In \cite{YY4}, the author  proved that the line solitary waves $Q_c(x-ct)$ of the Zakharov--Kuznetsov equation on $\RTL$ is orbitally stable and the asymptotically stable for $0<c\leq \frac{4}{5L^2}$ and is orbitally unstable for $c> \frac{4}{5L^2}$ by applying the argument in \cite{M M 3,C M P S}.
Pelinovsky \cite{DP} proved the asymptotic stability of the transversely modulated solitary waves of the Zakharov--Kuznetsov equation on $\RTL$ with exponentially wighted spaces.
For more results of the stability of line solitary waves to Kadomtsev and Petviashvili equation, we refer to the papers \cite{A P S,K P,TM2,TM3,R T 0,R T 1,R T 3,VEZ}, and the references therein.

In this paper, we construct center stable manifolds around unstable line solitary waves to the Zakharov--Kuznetsov equation on $\RTL$.
To introduce the main results, we define some notations.
The orbit of the solitary wave $Q_c$ is defined by 
\[S(c)=\{\tau_qQ_c: q \in \R\}\]
and the tubular neighborhood of the orbit of the solitary wave $Q_c$ is defined by
\[\mathcal{N}_{\delta,c}=\{u \in H^1(\RTL) : \inf_{q \in \R}\norm{u-\tau_qQ_c}_{H^1}<\delta\},\]
where $(\tau_qu)(x,y)=u(x-q,y)$.
By the global well-posedness result in \cite{M P}, we define $U(t)$ by the flow map of \eqref{ZKeq} at time $t$.
The following is the main theorem.
\begin{theorem}\label{thm-main}
Let $c^*>\frac{4}{5L}$. 
Then, there exists $C^1$ manifold $\mathcal{M}_{cs}(c^*)$ in $H^1(\RTL)$ containing the orbit $S(c^*)$ with the following properties:
\begin{enumerate}
  \setlength{\parskip}{0.1cm} 
  \setlength{\itemsep}{0.05cm} 
\renewcommand{\labelenumi}{\rm (\roman{enumi})}
\item The codimension of $\mathcal{M}_{cs}(c^*)$ in $H^1(\RTL)$ equals $2 \min\{n \in \Z: \frac{\sqrt{5c^*}L}{2}-1\}$ which is the total dimension of the eigenspaces of the linearized operator $\partial_x (-\Delta+c^*-2Q_{c^*})$ corresponding to eigenvalues with positive real part.
\item $\tau_qU(t)\mathcal{M}_{cs}(c^*) \subset \mathcal{M}_{cs}(c^*)$ for $q \in \R$ and $t\geq 0$.
\item $\mathcal{M}_{cs}(c^*)$ is normal at $Q_{c^*}$ to the corresponding to eigenvalues of $\partial_x (-\Delta+c^*-2Q_{c^*})$ with positive real part.
\item For any $\varepsilon>0$, there exists $\delta>0$ such that $U(t)(\mathcal{M}_{cs}(c^*)\cap \mathcal{N}_{\delta,c^*}) \subset \mathcal{M}_{cs}(c^*) \cap \mathcal{N}_{\varepsilon,c^*}$ for $t\geq 0$.
\item There is $\varepsilon_0>0$ such that for $u_0 \in \mathcal{N}_{\varepsilon_0,c^*}\setminus \mathcal{M}_{cs}(c^*)$ there exists $t_0>0$ satisfying $U(t_0)u_0 \notin \mathcal{N}_{\varepsilon_0,c^*}$.  
\end{enumerate}
\end{theorem}

\begin{remark}
In this paper, we only consider the solutions of \eqref{ZKeq} in \cite{M P} which are in the Bourgain spaces in local time.
If the unconditional uniqueness of the solutions to \eqref{ZKeq} in $C(\R, H^1(\RTL))$ is proved, we can obtain the center stable manifold without the restriction of the class of solutions.
\end{remark}

\begin{remark}
In \cite{YY5}, we obtain the  existence of the center stable manifold $\mathcal{M}_{cs}(c^*)$ for $c^* \in \{ c> 4/5L^2 : c \neq 4n^2/5L^2 \mbox{ for } n \in \Z\}$. 
Since the linearized operator $\partial_x (-\Delta+c^*-2Q_{c^*})$ with $c^*=4n^2/5L^2$ has extra eigenfunctions corresponding to $0$ eigenvalue for any positive integer $n$, it is difficult to obtain Theorem \ref{thm-main} by applying the argument in \cite{YY5} to the case $c^*=4n^2/5L^2$ directly.
\end{remark}

By developing the Hadamard method, Bates and Jones \cite{B J} proved a general theorem for the existence of invariant manifolds of equilibria for nonlinear partial differential equations.
In \cite{B J}, applying the general theorem and using a estimate of a Lyapunov function, Bates and Jones showed the existence of a Lipschitz center stable manifold of a stationary solution for nonlinear Klein--Gordon equation on $\R^n$ with the power nonlinearity $|u|^{p}u$ ($0<p<\frac{2}{n-2}$) under the radial symmetry restriction and the assumption which is the triviality of the null space of the linearized operator.
Nakanishi and Schlag \cite{N S 2} constructed center stable manifolds of ground states for nonlinear Klein--Gordon equation on $\R^n$ with the $H^1$-subcritical power nonlinearity without the radial symmetry restriction.
In \cite{B J}, Bates and Jones assume the Lipschitz continuity of the function $|u|^pu:H^1(\R^n) \to L^2(\R^n)$.
Using the Strichartz norm to the contraction of center stable manifold, Nakanishi and Schlag treated the $H^1$-subcritical power nonlinearity in \cite{N S 2}.
In the case of the non-radial symmetry, the linearized operator around ground state has non-trivial null space which comes from the translation symmetry and yields a the derivative loss in the Hadamard graph contraction argument.
To treat the derivative loss term due to the translation, in \cite{N S 2} Nakanishi and Schlag introduce the mobile distance for the construction of center stable manifolds.
By using the Strichartz estimate of the linear evolution around ground states, Schlag \cite{WS} constructed center stable manifolds in $W^{1,1} \cap W^{1,2}$ of ground states for the 3D cubic nonlinear Schr\"odinger equation and proved the asymptotic behavior of solutions on the center stable manifolds.
Applying the argument in \cite{WS}, Krieger and Schlag \cite{K S} constructed a center stable manifold of ground states for 1D nonlinear Schr\"odinger equation with $L^2$-supercritical nonlinearity.
The result \cite{WS} was improved by Beceanu \cite{MB} who constructed a center stable manifold of ground states for the 3D cubic nonlinear Schr\"odinger equation on the critical space $\dot{H}^{\frac{1}{2}}(\R^3)$.
By proving the trichotomy results which classifies initial datum near one soliton by the asymptotic behavior of solution, Martel, Merle, Nakanishi and Rapha\"el \cite{M M N R} constructed a center stable manifold of one soliton for the $L^2$-critical generalized KdV equation.
By using a smooth bundle coordinate system instead of a translational parametrization, Jin, Lin and Zeng \cite{J L Z 2} constructed a center stable manifold of one soliton for $L^2$-supercritical generalized KdV equation.
In \cite{J L Z 2}, applying smoothing estimate for solution of the Airy equation, Jin, Lin and Zeng treat the nonlinearity of generalized KdV equation with a loss of derivative.
In \cite{YY5}, the author constructed a center stable manifolds of line solitary waves $Q_{c^*}(x-c^*t)$ with  $c^* \in \{ c> 4/5L^2 : c \neq 4n^2/5L^2 \mbox{ for } n \in \Z\}$ for \eqref{ZKeq} by modifying the mobile distance in \cite{N S 2}.
Using the bi-linear estimate on Fourier restriction spaces in \cite{M P}, the author control a loss of derivative for the nonlinearity of \eqref{ZKeq}.
For more results of center stable manifolds of relative equilibria, we refer to the papers \cite{MB0,J L Z 1,K N S 0,K N S,N S}, and the references therein.

%In \cite{YY5}, we obtain the decomposition 
%\begin{align}\label{sec1-dec}
%u=P_+u+P_-u+P_0u+\gamma,
%\end{align}
%where $
Applying the argument in \cite{YY5} with a modulation of extra eigenfunctions corresponding to the $0$ eigenvalue, we obtain the existence of Lipschitz invariant manifolds of $Q_{c^*}(x-c^*t)$ with $c^* \in \{ 4n^2/5L^2 : n \in \Z, n > 1\}$.
For $c_0 \in \{ 4n^2/5L^2 : n \in \Z, n \geq 1\}$, the linearized operator of \eqref{ZKeq} around $Q_{c_0}(x-c_0t)$ has extra eigenfunctions corresponding to the $0$ eigenvalue and $Q_{c_0}(x-c_0t)$ is a bifurcation point of the stationary equation \eqref{eq-st} with the bifurcation parameter $c$.
We cannot control this degeneracy of the linearized operator by modulations for the translation symmetry of the equation \eqref{ZKeq}.
Therefore, we cannot use the coerciveness of the linearized operator around line solitary waves $Q_{c^*}(x-c^*t)$ on the invariant manifold to show the orbital stability of line solitary waves $Q_{c^*}(x-c^*t)$ with $c^* \in \{ 4n^2/5L^2 : n \in \Z, n > 1\}$ on the invariant manifold.
The orbital stability and the orbital instability of standing waves of nonlinear Schr\"odinger equation with the degeneracy of the linearized operator of the evolution equation was proved by Comech and Pelinovsky \cite{C P}, Maeda \cite{MM} and \cite{YY2}.
%To prove the orbital stability of standing waves with the degenerate kernel of the linearized operator,  
In the case that the linearized operator of the evolution equation has extra eigenfunctions corresponding to the $0$ eigenvalue, the positivity of the higher order term of the Taylor expansion for Lyapunov functions derives the orbital stability of standing waves.
In Section 3, to prove the orbital stability on the invariant manifold around $Q_{c^*}(x-c^*t)$, we apply the fourth order estimate of the Lyapunov function $S_{c}(u)$ for $Q_{4/5L^2}(x-4t/5L^2)$ in \cite{YY2,YY4}.
By the decomposition $u=P_+u+P_-u+P_0u+P_{\gamma}u$ around line solitary wave  $Q_{c^*}(x-c^*t)$ (see Proposition \ref{prop-linear} for the detail of the statement) and the estimate in \cite{N S 2,YY5}, we obtain that the order of unstable modes $P_+u$ on the invariant manifolds is controlled by the order of $u$ and the order of stable modes $P_-u$ on the invariant manifolds is controlled by the second order of $u$.
Since the dominant order of the positive term of the Lyapunov function $S_{c^*}(u)$ around the line solitary wave is $4$, by using the dominant order of the positive term of $S_{c^*}(u)$ we cannot control the error term which is the $L^2$ inner product of $P_+u$ and $\mathbb{L}_{c^*}P_-u$.
To obtain a sharper estimate of the error term, we construct invariant manifolds which is a Lipschitz graph function on the stable invariant space of the linearized operator with the H\"older exponent $\frac{3}{2}<\alpha<2$ at $Q_{c^*}(x-c^*t)$ (see the definition of $\mathcal{G}^+_{l_1,l_2,\alpha,\delta,\kappa}$ for the detail).
By the $\alpha$-H\"older continuity of the invariant manifolds at $Q_{c^*}(x-c^*t)$, we have that the order of unstable modes $P_+u$ on the invariant manifolds is $\alpha$ which yields the order of the error term is controlled by the $2+\alpha$-order of $u$.
Since $2+\alpha<4$, to show the order of the error term is controlled by the fourth order of $u$, we apply the bootstrap argument for the estimate of the order of the error term by using the Lyapunov function $S_{c^*}(u)$.

This paper is organized as following.
In Section 2, we introduce a spectral decomposition with respect to the linearized operator of \eqref{ZKeq} and the estimate of the difference between solutions of a localized equation of \eqref{ZKeq} by a mobile distance.
In Section 3, we construct invariant manifolds in $\mathcal{G}^+_{l_1,l_2,\alpha,\delta,\kappa}$ by applying the argument \cite{N S 2} and prove the orbital stability of line solitary wave on the invariant manifolds by using the fourth order estimate of Lyapunov function in \cite{YY4}.
In Section 4, we prove the $C^1$ regularity of the center stable manifolds.

\section{Linearized operator and localized equation}
In this section, we show the properties of the linearized operator around line solitary waves and introduce the localized equation around line solitary waves.
We define the linearized operator of the stationary equation
\[\mathbb{L}_c=-\Delta +c -2Q_c.\]
Let $c^*>4/5L^2$ and $n_0$ be the integer with $\frac{2(n_0-1)}{\sqrt{5c^*}}< L \leq \frac{2n_0}{\sqrt{5c^*}}$ and $n_0\geq 2$.
We define the set of critical speed by
\[\mbox{\rm CS}=\Bigl\{\frac{4n^2}{5L^2}; n \in \Z , n>1\Bigr\}.\]
The following proposition follows Proposition 3.1 in \cite{YY4} and Proposition 2.2 and Proposition 2.3 in \cite{YY5}.
\begin{proposition}\label{prop-linear}
Let $c^* > 4/5L^2$. The following holds.
\begin{enumerate}
  \setlength{\parskip}{0.1cm} 
  \setlength{\itemsep}{0.05cm} 
\renewcommand{\labelenumi}{\rm (\roman{enumi})}
\item Eigenvalues of $\partial_x\mathbb{L}_{c^*}$ with the positive real part are positive real number $\lambda_1,\lambda_2,\dots,\lambda_{n_0-1}$.
\item Eigenvalues of $\partial_x\mathbb{L}_{c^*}$ with the negative real part are negative real number \\
$-\lambda_1,-\lambda_2,\dots,-\lambda_{n_0-1}$.
\item There exist $f_1,f_2,\dots,f_{n_0-1} \in H^{\infty}(\R)$ such that for $k\in \{1,2,\dots,n_0-1\}$, 
\[F_k^{\pm,0}(x,y)=\pm f_k(\pm x)\cos \frac{ky}{L}, \quad F_k^{\pm,1}(x,y)=\pm f_k(\pm x)\sin \frac{ky}{L}\]
satisfy $\partial_x\mathbb{L}_{c^*}F_k^{\pm,j}=\pm\lambda_kF_k^{\pm,j}$ and $(F_k^{\pm,j},\mathbb{L}_{c^*}F_k^{\mp,j})_{L^2}=1$.
\item $\partial_x\mathbb{L}_{c^*}\partial_cQ_{c^*}=-\partial_xQ_{c^*}$ and 
\[\mbox{\rm Ker}(\partial_x\mathbb{L}_{c^*})=
\begin{cases}
\mbox{\rm Span}\{\partial_xQ_{c^*}\}, & c^* \notin \mbox{\rm CS},\\
\mbox{\rm Span}\{\partial_xQ_{c^*}, Q_{c^*}^{\frac{3}{2}}\cos \frac{n_0y}{L}, Q_{c^*}^{\frac{3}{2}}\sin \frac{n_0y}{L}\}, & c^*= \frac{4n_0^2}{5L^2} \in \mbox{\rm CS},
\end{cases}\]
where $\mbox{\rm Span}A$ is the linear subspace which is spanned by elements in the set $A$.
\item There exists $C>0$ such that for $u \in H^1(\RTL)$, 
\[P_{\gamma}u=u-\sum_{\substack{ j=0,1 \\ k=1,2, \dots, n_0-1}}(\Lambda_k^{+,j}F_k^{+,j}+\Lambda_k^{-,j}F_k^{-,j})-P_0u\]
satisfies 
\[\tbr{P_{\gamma}u,\mathbb{L}_{c^*}P_{\gamma}u}_{H^1,H^{-1}}\geq C \norm{P_{\gamma}u}_{H^1}^2,\]
where $\partial_cQ_{c^*}=\frac{\partial Q_c}{\partial c} |_{c=c^*}$,
\begin{align*}
&P_0u=\begin{cases}
\mu_1\partial_xQ_{c^*}+\mu_2\partial_cQ_{c^*}, & c^* \notin \mbox{\rm CS},\\
\mu_1\partial_xQ_{c^*}+\mu_2\partial_cQ_{c^*}+a_0 Q_{c^*}^{\frac{3}{2}}\cos \frac{n_0y}{L}+a_1 Q_{c^*}^{\frac{3}{2}}\sin \frac{n_0y}{L}, & c^*= \frac{4n_0^2}{5L^2} \in \mbox{\rm CS}, 
\end{cases}\\
&\Lambda_k^{\pm,j}=(u,\mathbb{L}_{c^*}F_k^{\mp,j})_{L^2}, \ \mu_1=\frac{(u,\partial_xQ_{c^*})_{L^2}}{\norm{\partial_xQ_{c^*}}_{L^2}^2}, \ \mu_2=\frac{(u,Q_{c^*})_{L^2}}{(\partial_cQ_{c^*},Q_{c^*})_{L^2}}, \\
& a_0=\frac{(u,Q_{c^*}^{\frac{3}{2}}\cos \frac{n_0 y}{L})_{L^2}}{\norm{Q_{c^*}^{\frac{3}{2}}\cos \frac{n_0 y}{L}}_{L^2}^2}, \quad a_1=\frac{(u,Q_{c^*}^{\frac{3}{2}}\sin \frac{n_0 y}{L})_{L^2}}{\norm{Q_{c^*}^{\frac{3}{2}}\sin \frac{n_0 y}{L}}_{L^2}^2}.
\end{align*}
\end{enumerate}
\end{proposition}
We define the spectral projection
\[P_{\pm}u=\sum_{\substack{ j=0,1 \\ k=1,2, \dots, n_0-1}}\Lambda_k^{\pm,j}F_k^{\pm,j}, \quad P_1u=\mu_1\partial_xQ_{c^*}, \quad P_2u=\mu_2\partial_cQ_{c^*},\]
\[P_au=a_0Q_{c^*}^{\frac{3}{2}}\cos \frac{n_0 y}{L}+a_1Q_{c^*}^{\frac{3}{2}}\sin \frac{n_0 y}{L}, \quad P_du=u-P_{\gamma}u\]
and the norm
\[\norm{u}_{E_\kappa}=
\Bigl(\sum_{\substack{ j=0,1 \\ k=1,2, \dots, n_0-1}}( (\Lambda_k^{+,j})^2+(\Lambda_k^{-,j})^2) + \kappa^2 \mu_1^2+\mu_2^2+ \kappa^2 1_{\mbox{\rm CS}}(c^*)(a_0^2+a_1^2)+ \tbr{\gamma,\mathbb{L}_{c^*}\gamma}_{H^1,H^{-1}} \Bigr)^{\frac{1}{2}}, \]
where $\kappa>0$, $\Lambda_k^{\pm,j}, \mu_1,\mu_2, a_0$ and $a_1$ are given in Proposition \ref{prop-linear} and $1_A$ is the indicator function of $A$.

Let $u$ be a solution to the equation \eqref{ZKeq}.
Then, $v(t)=\tau_{-\rho(t)}u(t)-Q_{c(t)}$ satisfies
\begin{equation}\label{v-eq}
v_t=\partial_x\mathbb{L}_{c^*}v+(\dot{\rho}-c)\partial_xQ_{c^*}-\dot{c}\partial_cQ_{c^*}+N(v,c,\rho),
\end{equation}
where 
\[N(v,c,\rho)=\partial_x[-v^2+(\dot{\rho}-c^*)v+2(Q_{c^*}-Q_c)v+(\dot{\rho}-c)(Q_c-Q_{c^*})]-\dot{c}\partial_c(Q_c-Q_{c^*}).\]
In the following lemma, we choose the modulation parameters $\rho$ and $c$ satisfying the orthogonality condition
\begin{equation}\label{orth-cond-1}
(v,\partial_xQ_{c^*})_{L^2}=(v,\partial_cQ_{c^*})_{L^2}=0.
\end{equation}
\begin{lemma}\label{lem-orth}
There exist $\delta_0, C_{\delta_0}>0$ and smooth maps $\rho:\mathcal{N}_{\delta_0,c^*} \to \R$ and $c:\mathcal{N}_{\delta_0,c^*} \to (0,\infty)$ such that for $u \in \mathcal{N}_{\delta_0,c^*}$, $v=\tau_{-\rho(u)}u-Q_{c(u)}$ satisfies orthogonality condition $\eqref{orth-cond-1}$ and 
\begin{equation*}
\norm{v}_{H^1}+|c(u)-c^*|< C_{\delta_0} \inf_{q\in \R}\norm{u-\tau_qQ_{c^*}}_{H^1},
\end{equation*}
where $\mathcal{N}_{\delta,c^*}=\{u \in H^1(\RTL); \inf_{q \in \R}\norm{u-\tau_qQ_{c^*}}_{H^1}<\delta\}$.
\end{lemma}
\proof
We define
\[G(u,c,\rho)=
\begin{pmatrix}
(\tau_{-\rho}u-Q_c,\partial_x Q_{c^*})_{L^2}\\
(\tau_{-\rho}u-Q_c,Q_{c^*})_{L^2}
\end{pmatrix}.\]
Then, $G(Q_{c^*},c^*,0)={}^t(0,0)$ and $\frac{\partial G}{\partial c \partial \rho}(Q_{c^*},c^*,0)$ is regular.
By the implicit function theorem, we obtain the conclusion.

\qed

For solution $(v(t),c(t),\rho(t))$ to the equation \eqref{v-eq} satisfying \eqref{orth-cond-1}, $(c(t),\rho(t))$ satisfies
\begin{equation}\label{c-rho-eq}
\begin{pmatrix}
\dot{\rho}-c\\
\dot{c}
\end{pmatrix}
=
\begin{pmatrix}
\norm{\partial_xQ_{c^*}}_{L^2}^{-2}(v,\mathbb{L}_{c^*}\partial_x^2Q_{c^*})_{L^2}\\
0
\end{pmatrix}
+\bm{N}(v,c),
\end{equation}
where
\begin{align*}
\bm{N}(v,c)
=&
\begin{pmatrix}
(\partial_xQ_{c}+\partial_x v, \partial_xQ_{c^*})_{L^2} & 0\\
0 & -(\partial_cQ_{c},Q_{c^*})_{L^2}
\end{pmatrix}^{-1}\\
 & \times 
\begin{pmatrix}
(\partial_x\mathbb{L}_{c^*}v+(c-c^*)\partial_xv-\partial_xv^2+2\partial_x((Q_{c^*}-Q_c)v),\partial_xQ_{c^*})_{L^2}\\
(-\partial_xv^2 +2\partial_x((Q_{c^*}-Q_c)v),Q_{c^*})_{L^2}
\end{pmatrix}\\
& - 
\begin{pmatrix}
\norm{\partial_xQ_{c^*}}_{L^2}^{-2}(v,\mathbb{L}_{c^*}\partial_x^2Q_{c^*})_{L^2} \\
0
\end{pmatrix}\\
=&O(\norm{v}_{L^2}^2+\norm{v}_{L^2}|c-c^*|) \mbox{ as } \norm{v}_{L^2}+|c-c^*| \to 0.
\end{align*}

On the tubular neighborhood $\mathcal{N}_{\delta_0,c^*}$, $u=\tau_{\rho}(v+Q_c)$ is a solution to the equation \eqref{ZKeq} satisfying the orthogonality condition \eqref{orth-cond-1} if and only if $v=\tau_{-\rho}u-Q_c$ is solution to the equation \eqref{v-eq} with $(c,\rho)$ satisfying \eqref{c-rho-eq} and $(v(0),\partial_xQ_{c^*})_{L^2}=(v(0),Q_{c^*})_{L^2}=0$.

Let $\chi \in C_0^{\infty}(\R)$ be a smooth function with
\begin{equation*}
\chi(r)=
\begin{cases}
1, & |r|\leq 1 \\
0, & |r| \geq 2
\end{cases}
,\quad 
0\leq \chi \leq 1
\end{equation*}
and
\[\chi_{\delta}=\chi_{\delta}(v,c)=\chi(\delta^{-2}(\norm{v}_{H^1}^2+|c-c^*|^2)).\]
We define the localized system of the system \eqref{v-eq} and \eqref{c-rho-eq} as
\begin{equation}
v_t=\partial_x\mathbb{L}_{c^*}v+(\dot{\rho}-c)\partial_xQ_{c^*}-\dot{c}\partial_cQ_{c^*}+\chi_{\delta}(v,c)N(v,c,\rho), \label{Leq-1}
\end{equation}
\begin{equation}
\begin{pmatrix}
\dot{\rho}-c\\
\dot{c}
\end{pmatrix}
=
\begin{pmatrix}
\norm{\partial_xQ_{c^*}}_{L^2}^{-2}(v,\mathbb{L}_{c^*}\partial_x^2Q_{c^*})_{L^2}\\
0
\end{pmatrix}
+ \bm{N}_{\delta}(v,c), \label{Leq-2}
\end{equation}
where
\begin{align*}
&\bm{N}_{\delta}(v,c)\\
=&
\begin{pmatrix}
\norm{\partial_xQ_{c^*}}_{L^2}^2+\chi_{\delta}(\partial_x(v+Q_c-Q_{c^*}),\partial_xQ_{c^*})_{L^2} & 0 \\
0 & -(\partial_cQ_{c^*}+\chi_{\delta}\partial_c(Q_c-Q_{c^*}),Q_{c^*})_{L^2}
\end{pmatrix}^{-1}\\
& \times 
\begin{pmatrix}
\chi_{\delta}(\partial_x((c-c^*)v-v^2+2(Q_{c^*}-Q_c)v),\partial_xQ_{c^*})_{L^2}-(v,\mathbb{L}_{c^*}\partial_x^2Q_{c^*})_{L^2}\\
-\chi_{\delta}(\partial_x(v^2-2(Q_{c^*}-Q_c)v),Q_{c^*})_{L^2}
\end{pmatrix}\\
&-
\begin{pmatrix}
\norm{\partial_xQ_{c^*}}_{L^2}^{-2}(v,\mathbb{L}_{c^*}\partial_x^2Q_{c^*})_{L^2}\\
0
\end{pmatrix}\\
=&\chi_{\delta}O(\norm{v}_{L^2}^2+\norm{v}_{L^2}|c-c^*|) \mbox{ as } \norm{v}_{L^2}+|c-c^*| \to 0.
\end{align*}

To obtain solutions to the system \eqref{Leq-1} and \eqref{Leq-2}, we solve the system
\begin{align}
&w_t=-\partial_x\Delta w-2\partial_x((\tau_{\rho_*}Q_{c^*})w)+(\dot{\rho}-c)\tau_{\rho_*}\partial_xQ_{c^*}-\dot{c}\tau_{\rho_*}\partial_cQ_{c^*}+\chi_{\delta}(w,c)\tilde{N}(w,c,\rho), \label{w-eq-1}\\
&
\begin{pmatrix}
\dot{\rho}-c\\
\dot{c}
\end{pmatrix}
=
\begin{pmatrix}
\norm{\partial_xQ_{c^*}}_{L^2}^{-2}(w,\tau_{\rho_*}(\mathbb{L}_{c^*}\partial_x^2Q_{c^*}))_{L^2}\\
0
\end{pmatrix}
+\tilde{\bm{N}}_{\delta}(w,c,\rho), \label{w-eq-2}
\end{align}
where
\[\tilde{N}(w,c,\rho)=\partial_x[-w^2+2w\tau_{\rho_*}(Q_{c^*}-Q_c)+(\dot{\rho}-c)\tau_{\rho_*}(Q_c-Q_{c^*})]-\dot{c}\tau_{\rho_*}\partial_c(Q_c-Q_{c^*}),\]
$\tilde{\bm{N}}_{\delta}(w,c,\rho)=\bm{N}_{\delta}(\tau_{-\rho_*}w,c)$ and
\[\rho_*(w,c,\rho,t)=\rho_*(t)=c^*t+\int_0^t \chi_{\delta}(w(s),c(s)-c^*)(\dot{\rho}(s)-c^*)\, ds.\]
Then, we have $v=\tau_{\rho_*}w$ and the system \eqref{w-eq-1}--\eqref{w-eq-2} has no the advection term $\partial_x w$.
To solve the system \eqref{w-eq-1}--\eqref{w-eq-2}, we define the Bourgain space $X^{s,b}$ related to the linear part of \eqref{ZKeq} as the completion of the Schwartz space under the norm
\[\norm{u}_{X^{s,b}}=\Bigl( \int_{\R^2} \sum_{L\eta \in \Z} \tbr{\tau-\xi(\xi^2+\eta^2)}^{2b}\tbr{\sqrt{3\xi^2+\eta^2}}^{2s}|\tilde{u}(\tau,\xi,\eta)|^2d\tau d\xi\Bigr)^{\frac{1}{2}},\]
where $\tbr{x}=1+|x|$ and $\tilde{u}$ is the space-time Fourier transform of $u$.
Moreover, for $T_1<T_2$ we define the localized space $X_{T_1,T_2}^{s,b}$ of $X^{s,b}$ by the norm
\[\norm{u}_{X_{T_1,T_2}^{s,b}}=\inf\{\norm{v}_{X^{s,b}};v \in X^{s,b}, v(t)=u(t) \mbox{ for } t \in [T_1,T_2]\}.\]
Particularly, we denote $X_{-T,T}^{s,b}$ by $X_T^{s,b}$.
The following theorem shows the global well-posedness of the systems \eqref{Leq-1}--\eqref{Leq-2} and \eqref{w-eq-1}--\eqref{w-eq-2}.
\begin{theorem}\label{thm-w-gwp}
Let $0<\kappa<1$. The system $\eqref{Leq-1}$--$\eqref{Leq-2}$ is globally well-posed in $H^1(\RTL)\times (0,\infty) \times \R$.
Precisely, there exists $b>1/2$ such that for every $(v_0,c_0,\rho_0) \in H^1(\RTL)\times (0,\infty) \times \R$ and $T>0$ there exists a unique solution $(w,c,\rho)$ of the system $\eqref{w-eq-1}$--$\eqref{w-eq-2}$ such that $(w(0),c(0),\rho(0))=(v_0,c_0,\rho_0)$, $(w,\dot{c},\dot{\rho}-c) \in X_T^{1,b} \times L^2(-T,T) \times L^2(-T,T)$ and $(\tau_{-\rho_*}w,c,\rho)$ is a solution to the system $\eqref{Leq-1}$--$\eqref{Leq-2}$ with initial data $(v(0),c(0),\rho(0))$.
Moreover, the flow map of the system $\eqref{w-eq-1}$--$\eqref{w-eq-2}$ is Lipschitz continuous on bounded sets of $H^1(\RTL) \times (0,\infty) \times \R$ and there exists $T^*,C_\kappa>0$ such that for any $0<\delta<1$ and initial data $(v(0),c(0),\rho(0))$ satisfies
\begin{align}
\sup_{|t|\leq T^*} \norm{v(t)}_{E_\kappa} +\norm{\dot{c}}_{L^2(-T^*,T^*)}+\norm{\dot{\rho}-c}_{L^2(-T^*,T^*)} \leq C_\kappa \norm{v(0)}_{E_\kappa},\notag \\
\sup_{|t|\leq T^*} \norm{P_d(v(t)-e^{t\mathcal{A}}v(0))}_{E_\kappa} \leq C_\kappa \min\{\norm{v(0)}_{E_\kappa}(\norm{v(0)}_{E_\kappa}+|c(0)-c^*|),\delta^2 \},\notag \\
\sup_{|t|\leq T^*} |\norm{P_{\gamma}v(t)}_{E_\kappa}^2-\norm{P_{\gamma}v(0)}_{E_\kappa}^2| \leq C_\kappa \min\{ \norm{v(0)}_{E_\kappa}^2(\norm{v(0)}_{E_\kappa}+|c(0)-c^*|), \delta^3\} \label{eq-w-gwp-1},
\end{align}
where the constants $C_\kappa$ and $T^*$ do not depend on $\delta$ and $(v(0),c(0),\rho(0))$ and
\[\mathcal{A}u=\partial_x\mathbb{L}_{c^*}u+\frac{(u,\mathbb{L}_{c^*}\partial_x^2Q_{c^*})_{L^2}}{\norm{\partial_xQ_{c^*}}_{L^2}^{2}}\partial_xQ_{c^*}.\]
\end{theorem}
To show the global well-posedness of the system $\eqref{Leq-1}$--$\eqref{Leq-2}$ in Theorem \ref{thm-w-gwp}, we use the estimates 
\begin{equation}\label{eq-w-gwp-2}
\norm{\partial_x(Qu)}_{X^{1,0}} \lesssim \Bigl(\norm{\partial_xQ}_{L^{\infty}_t W_{x,y}^{1,\infty}}+\sum_{|\alpha|\leq 1}\norm{\partial^\alpha Q}_{L_x^2L_{yt}^{\infty}}\Bigr)\norm{u}_{X^{1,b}}
\end{equation}
and
\begin{equation}\label{eq-w-gwp-3}
\norm{\partial_x(uv)}_{X^{1,-\frac{1}{2}+2\varepsilon}}\lesssim \norm{u}_{X^{1,\frac{1}{2}+\varepsilon}}\norm{v}_{X^{1,\frac{1}{2}+\varepsilon}}
\end{equation}
by Molinet--Saut--Tzvetkov \cite{M S T} and Molinet--Pilod \cite{M P}, where $b>1/2$, $0<\varepsilon \ll 1$, $\alpha=(\alpha_1,\alpha_2)$ and $\partial^\alpha=\partial_x^{\alpha_1}\partial_y^{\alpha_2}$.
Moreover, combining  the estimates \eqref{eq-w-gwp-2}, \eqref{eq-w-gwp-3} and 
\begin{equation}\label{eq-w-gwp-4}
\left| \int_0^T \int_{\RTL} (-\Delta u)v dtdxdy \right|\lesssim \norm{\tau_{\rho}u}_{X^{1,b}}\norm{\tau_{\rho}v}_{X^{1,-\beta}}
\end{equation}
for $\tau_{\rho}u \in X^{1,b}, \tau_{\rho}v \in X^{1,-\beta},\rho \in L^{\infty}_t,p>1, b >\frac{1}{2}$ and $0\leq \beta <\frac{1}{2}$ with $(b-\beta)p>p-1$, for any solution $(v,c,\rho)$ to the system $\eqref{Leq-1}$--$\eqref{Leq-2}$ and solution $(w,c,\rho)$ of the system $\eqref{w-eq-1}$--$\eqref{w-eq-2}$ we obtain
\begin{align}\label{est-cont-sp}
\sup_{|t|\leq T^*} |\norm{P_{\gamma}v(t)}_{E_\kappa}^2-\norm{P_{\gamma}v(0)}_{E_\kappa}^2| \leq& C\norm{w}_{X^{1,b}}^2(\norm{w}_{X^{1,b}} + |c(0)-c^*|)\notag \\
\lesssim& \norm{v(0)}_{H^1}^2(\norm{v(0)}_{H^1}+|c(0)-c^*|)
\end{align}
and the estimate \eqref{eq-w-gwp-1}, where the constant $C$ does not depend on $\kappa$.
Since the proof of Theorem \ref{thm-w-gwp} follows the proof of Theorem 3.4 in \cite{YY5}, we omit the detail of the proof.

To construction a invariant manifolds from the flow map of the system \eqref{Leq-1}--\eqref{Leq-2}, we show estimates of the difference between solutions to the system \eqref{Leq-1}--\eqref{Leq-2} and solutions to the linearized equation of \eqref{Leq-1}--\eqref{Leq-2}.
Theorem \ref{thm-w-gwp} yields the Lipschitz continuity of the flow map of the system \eqref{w-eq-1}--\eqref{w-eq-2} by the energy norm.
However, Theorem \ref{thm-w-gwp} does not imply the Lipschitz continuity of the flow map of the system \eqref{Leq-1}--\eqref{Leq-2} by the energy norm.
To show the estimates, we define the mobile distance which was introduce in \cite{N S 2}.
Let $C_2$ be a large real constant and $\phi$ be the smooth positive non-deceasing function with
\[\phi(r)=\begin{cases}
1, & r \leq C_2,\\
r, & r \geq 2C_2.
\end{cases}\]
We define $\phi_{\delta}$ by 
\[\phi_{\delta}(u)=\phi(\delta^{-1}\norm{P_{\gamma}u}_{E_\kappa}).\]
\begin{definition}
Let $\delta, \kappa>0$.
We define the mobile distance $\mathfrak{m}_{\delta,\kappa}:(H^1(\RTL)\times (0,\infty))^2 \to [0,\infty)$ by
\begin{align*}
\mathfrak{m}_{\delta,\kappa}(\bm{v}_0,\bm{v}_1)=&\Bigl[ \norm{P_d(v_0-v_1)}_{E_\kappa}^2+ \inf_{q \in \R, j=0,1}(\norm{P_{\gamma}v_j-\tau_qP_{\gamma}v_{1-j}}_{E_\kappa}^2+\delta |q|^2\phi_{\delta}(v_{1-j})^2)\\
&+|\log c_0 -\log c_1|^2\Bigl]^{\frac{1}{2}}
\end{align*}
for $\bm{v}_0=(v_0,c_0),\bm{v}_1=(v_1,c_1) \in H^1(\RTL) \times (0,\infty)$.
\end{definition}
In the following lemma, we recall the mobile distance $\mathfrak{m}_{\delta,\kappa}$ is a quasi-distance on $H^1(\RTL)\times (0,\infty)$. 
\begin{lemma}\label{lem-dist}
Let $0<\delta,\kappa<1$. Then, $\mathfrak{m}_{\delta,\kappa}$ satisfies the following.
\begin{enumerate}
  \setlength{\parskip}{0.1cm} 
  \setlength{\itemsep}{0.05cm} 
\renewcommand{\labelenumi}{\rm (\roman{enumi})}
\item $\mathfrak{m}_{\delta,\kappa}(\bm{v}_0,\bm{v}_1)=\mathfrak{m}_{\delta,\kappa}(\bm{v}_1,\bm{v}_0)\geq 0$, where the equality holds iff $\bm{v}_0=\bm{v}_1$.
\item $\mathfrak{m}_{\delta,\kappa}(\bm{v}_0,\bm{v}_1)\leq C(\mathfrak{m}_{\delta,\kappa}(\bm{v}_0,\bm{v}_2)+\mathfrak{m}_{\delta,\kappa}(\bm{v}_2,\bm{v}_1))$, for some absolute constant $C>0$ which does not depend on $\delta$.
\item If $\mathfrak{m}_{\delta,\kappa}(\bm{v}_n,\bm{v}_m) \to 0$ $(n,m \to \infty)$, then $\{\bm{v}_n\}_n$ converges in $H^1(\RTL)\times(0,\infty)$.
\item For $\bm{v}_0=(v_0,c_0)$, $\bm{v}_1=(v_1,c_1) \in H^1(\RTL)\times (0,\infty)$
\begin{align*}
|\norm{v_0}_{H^1}-\norm{v_1}_{H^1}|+\norm{v_0-v_1}_{L^2}+|\log c_0 -\log c_1| \lesssim & \mathfrak{m}_{\delta,\kappa}(\bm{v}_0,\bm{v}_1) \\
\lesssim & \norm{v_0-v_1}_{H^1} + |\log c_0-\log c_1|,
\end{align*}
where the implicit constants do not depend on $\delta$.
\end{enumerate}
\end{lemma}
Since 
\[|q|\norm{\nabla P_{\gamma} v_{1-j}}_{L^2} \lesssim \delta^{1/2}|q|\phi_{\delta}(v_{1-j}),\]
Lemma \ref{lem-dist} follows Proposition 2.2 in \cite{N S 2}.

In the following lemma, we show the estimates between solutions to the system \eqref{Leq-1}--\eqref{Leq-2} and solutions to the linearized equation of \eqref{Leq-1}--\eqref{Leq-2} on the quasi-metric space $(H^1(\RTL)\times(0,\infty),\mathfrak{m}_{\delta,\kappa})$.
\begin{lemma}\label{lem-est-mobile}
Let $0 <\kappa <1$.
There exists $T^*, \delta^*, C_\kappa>0$ such that for any $0<\delta<\delta^*$ and solutions $(\bm{v}_j,\rho_j)=(v_j,c_j,\rho_j)$ to the system $\eqref{Leq-1}$--$\eqref{Leq-2}$ given in Theorem $\ref{thm-w-gwp}$, we have
\begin{align*}
&\sup_{|t|\leq T^*} \mathfrak{m}_{\delta,\kappa}(\bm{v}_0(t),\bm{v}_1(t)) \leq C_\kappa \mathfrak{m}_{\delta,\kappa}(\bm{v}_0(0),\bm{v}_1(0)),\\
&\sup_{|t|\leq T^*} \Bigl( \norm{P_d(v_0(t)-v_1(t)-e^{t\mathcal{A}}(v_0(0)-v_1(0)))}_{E_\kappa}^2 +|I(v_0(t),v_1(t))-I(v_0(0),v_1(0))|\Bigr) \\
\leq  &  C_\kappa \delta^{1/2} \mathfrak{m}_{\delta,\kappa}(\bm{v}_0(0),\bm{v}_1(0))^2,
\end{align*}
where  
\[I(v_0,v_1)=\inf_{q\in \R, j=0,1}(\norm{P_{\gamma}v_j-\tau_q P_{\gamma}v_{1-j}}_{E_\kappa}^2 + \delta |q|^2 \phi_{\delta}(v_{1-j})^2).\]
\end{lemma}
In the case $c^* \notin \{\frac{4n^2}{5L^2}; n \in \Z, n>1\}$, Lemma \ref{lem-est-mobile} is same as Lemma 3.8 in \cite{YY5}.
Since the estimate of $P_a(v_0(t)-v_1(t)-e^{t\mathcal{A}}(v_0(0)-v_1(0)))$ follows the local well-posedness in Theorem \ref{thm-w-gwp}, for $c^* \in \{\frac{4n^2}{5L^2}; n \in \Z, n>1\}$ the proof of Lemma \ref{lem-est-mobile} also follows the arguments in the proof of Lemma 3.2 in \cite{N S 2} and the proof of Lemma 3.8 in \cite{YY5}.
Therefore, we omit the proof of Lemma \ref{lem-est-mobile}.

\section{Construction of the center stable manifolds}
In this section, we construct the center stable manifolds by applying the argument in \cite{N S 2}.

Let $\mathcal{H}=:H^1(\RTL)\times (0,\infty)$.
For $l_1, l_2,\alpha,\delta,\kappa>0$, we define the set of the graph function on the stable invariant space of the linearized operator $\mathcal{A}$ with the H\"older exponent $\alpha$ at $(0,c^*)$ by
\begin{align*}
\mathscr{G}_{l_1,l_2,\alpha,\delta,\kappa}^+=\{ &G:\mathcal{H}  \to P_+H^1(\RTL);  G =G\circ P_{\leq 0}, G(0,c^*)=0, \\
&\norm{G(\bm{v})}_{E_{\kappa}}\leq l_1 \norm{\bm{v}}_{E_{\kappa}}^{\alpha} \mbox{ for } \bm{v} \in \mathcal{H},\\
& \norm{G(\bm{v}_0)-G(\bm{v}_1)}_{E_{\kappa}}\leq l_2 \mathfrak{m}_{\delta,\kappa}(\bm{v}_0,\bm{v}_1) \mbox{ for } \bm{v}_0,\bm{v}_1 \in \mathcal{H}\},
\end{align*}
where $P_{\leq 0}(v,c)=((I-P_+)v,c)$ and
\[\norm{(v,c)}_{E_\kappa}=(\norm{v}_{E_\kappa}^2+|\log c-\log c^*|^2)^{1/2}.\]
We define the graph $\gbr{G}$ of $G \in \mathscr{G}_{l,\delta,\kappa}^+$ as 
\[\{(v,c) \in \mathcal{H}; P_+v=G(v,c)\}.\]

\begin{lemma}\label{lem-linear-eq}
For any solution $u(t)$ to $u_t=\mathcal{A}u$ satisfies
\begin{equation}\label{lin-eq-1}
(u(t),\partial_xQ_{c^*})_{L^2}=(u(0),\partial_xQ_{c^*})_{L^2}, \quad
% \end{align}
% and
% \begin{align}\label{lin-eq-2}
(u(t),Q_{c^*})_{L^2}=(u(0),Q_{c^*})_{L^2}.
\end{equation}
% where
% \[\mathcal{A}u=\partial_x\mathbb{L}_{c^*}u+\frac{(u,\mathcal{L}_{c^*}\partial_x^2Q_{c^*})_{L^2}}{\norm{\partial_xQ_{c^*}}_{L^2}^2}\partial_xQ_{c^*}.\]
\end{lemma}
\proof
The equalities \eqref{lin-eq-1} follow
\[(\mathcal{A}u,\partial_xQ_{c^*})_{L^2}=(\mathcal{A}u,Q_{c^*})_{L^2}=0.\]
\qed

In the following lemma, we prove the upper estimate of the unstable eigen mode.
\begin{lemma}\label{lem-lin-est}
There exist $T^*>0$, $C_{\kappa}>0$ and $C_L>0$ such that if $l_2,\delta,\kappa>0$ satisfy 
\begin{equation}\label{ass-small}
l_2+C_\kappa \delta+\kappa \ll 1 \mbox{ and } C_\kappa l_2^{-1}\delta^{\frac{1}{4}} \ll  1,
\end{equation}
then for any solutions $(v_j,c_j,\rho_j)$ to the system $\eqref{Leq-1}$--$\eqref{Leq-2}$ $(j=0,1)$ satisfying 
\begin{equation}\label{ass-lin}
\norm{P_+(v_0(0)-v_1(0))}_{E_\kappa}\leq l_2 \mathfrak{m}_{\delta,\kappa}(\bm{v}_0(0),\bm{v}_1(0))
\end{equation}
one has
\begin{equation}\label{lin-con}
\norm{P_+(v_0(t)-v_1(t))}_{E_\kappa}\leq 
\begin{cases}
C_Ll_2\mathfrak{m}_{\delta,\kappa}(\bm{v}_0(t),\bm{v}_1(t)), & |t|\leq T^*, \\
l_2 \mathfrak{m}_{\delta,\kappa}(\bm{v}_0(t),\bm{v}_1(t)), & -T^*\leq t \leq -\frac{T^*}{2}.
\end{cases}
\end{equation}
\end{lemma}
\proof
In the case $c^* \notin \{\frac{4n^2}{5L^2}; n \in \Z, n>1\}$, the conclusion follows Lemma 4.1 in \cite{YY5}.
Therefore, we consider the case $c^* =\frac{4n_0^2}{5L^2}$ for some $n_0 \in  \{n \in \Z;n >1\}$.
% By Lemma \ref{lem-linear-eq}, we have
% \[(e^{t\mathcal{A}}(v_0(0)-v_1(0)),\partial_xQ_{c^*})_{L^2}=(v_0(0)-v_1(0),\partial_xQ_{c^*})_{L^2}\]
% and
% \[(e^{t\mathcal{A}}(v_0(0)-v_1(0)),Q_{c^*})_{L^2}=(v_0(0)-v_1(0),Q_{c^*})_{L^2}.\]
Let 
\begin{equation}\label{def-k}
k_*=\min_{k=1,2,\dots,n_0-1}\lambda_k, \quad k^*=\max_{k=1,2,\dots,n_0-1}\lambda_k.
\end{equation}
By Lemma \ref{lem-est-mobile} and \eqref{ass-lin}, we have
\begin{equation}\label{lin-eq-3}
\norm{P_+(v_0(t)-v_1(t))}_{E_{\kappa}}\leq (\max\{e^{k^*t},e^{k_*t}\}l_2+C_{\kappa}\delta^{\frac{1}{4}}) \mathfrak{m}_{\delta,\kappa}(\bm{v}_0(0),\bm{v}_1(0)),
\end{equation}
\begin{equation}\label{lin-eq-3-1}
\bigl||\log c_0(t)-\log c_1(t)|^2-|\log c_0(0) -\log c_1(0)|^2\bigr| \lesssim \delta \mathfrak{m}_{\delta,\kappa}(\bm{v}_0(0),\bm{v}_1(0))^2,
\end{equation}
and
\begin{align}\label{lin-eq-4}
&\mathfrak{m}_{\delta,\kappa}(\bm{v}_0(t),\bm{v}_1(t))^2-|\log c_0(t)-\log c_1(t)|^2 \notag \\
\geq& \norm{P_d(e^{t\mathcal{A}}(v_0(0)-v_1(0)))}_{E_\kappa}^2+I(v_0(0),v_1(0))-C_{\kappa}^2 \delta^{\frac{1}{2}}\mathfrak{m}_{\delta,\kappa}(\bm{v}_0(0),\bm{v}_1(0))^2
\end{align}
for sufficiently small $|t|$.
Since $ (P_1+P_a)\mathbb{L}_{c^*}=0$, for $|t|<1$ there exists $C>0$ such that
\begin{equation}\label{lin-eq-5}
\norm{P_ae^{t\mathcal{A}}(v_0(0)-v_1(0))}_{E_\kappa}\leq \norm{P_a (v_0(0)-v_1(0))}_{E_\kappa}+ C \kappa \mathfrak{m}_{\delta,\kappa}(\bm{v}_0(0),\bm{v}_1(0)).
\end{equation}
Combining the estimates \eqref{lin-eq-3-1}--\eqref{lin-eq-5}, the assumption \eqref{ass-lin} and Lemma \ref{lem-linear-eq}, we have there exist $C,C_\kappa,T>0$ such that 
\begin{align}\label{lin-eq-6}
\mathfrak{m}_{\delta,\kappa}(\bm{v}_0(t),\bm{v}_1(t))^2 \geq 
\begin{cases}
(1-l_2^2+e^{2k^*t}l_2^2-C\kappa-C_{\kappa}\delta^{\frac{1}{2}})\mathfrak{m}_{\delta,\kappa}(\bm{v}_0(0),\bm{v}_1(0))^2, & -T\leq t \leq 0\\
(e^{-2k^*|t|}-C\kappa-C_\kappa\delta^{\frac{1}{2}})\mathfrak{m}_{\delta,\kappa}(\bm{v}_0(0),\bm{v}_1(0))^2, & |t|\leq T.
\end{cases}
\end{align}
The inequalities \eqref{lin-eq-3} and \eqref{lin-eq-6} yield the inequalities \eqref{lin-con} for sufficiently small $\kappa$ and $\delta$.

\qed

In the following lemma, we prove the upper estimate of the unstable mode around $0$.
\begin{lemma}\label{lem-lip-est}
%Let $c^* =\frac{4n_0^2}{5L^2}$ for some $n_0 \in  \{n \in \Z;n >1\}$
For $1<\alpha<2$, there exist $T^*>0$, $C_{\kappa}>0$ and $C_\alpha>0$ such that if $l_1,\delta,\kappa>0$ satisfy 
\begin{equation}\label{ass-small-2}
l_1+C_\kappa \delta+\kappa \ll 1 \mbox{ and } C_\kappa l_1^{-1}\delta^{\frac{1}{4}} + C_{\kappa}l_1^{-1}\delta^{2-\alpha}\ll  1,
\end{equation}
then for any solutions $(v,c,\rho)$ to the system $\eqref{Leq-1}$--$\eqref{Leq-2}$ satisfying 
\begin{equation}\label{ass-lin-2}
\norm{P_+v(0)}_{E_\kappa}\leq l_1 \min\{ \norm{\bm{v}(0)}_{E_{\kappa}}^\alpha, \norm{\bm{v}(0)}_{E_{\kappa}}\}
\end{equation}
one has
\begin{equation}\label{lin-con-2}
\norm{P_+v(t)}_{E_\kappa}\leq 
\begin{cases}
C_\alpha l_1\norm{\bm{v}(t)}_{E_{\kappa}}^\alpha, & |t|\leq T^*, \\
l_1 \norm{\bm{v}(t)}_{E_{\kappa}}^\alpha, & -T^*\leq t \leq -\frac{T^*}{2},
\end{cases}
\end{equation}
where $\bm{v}(t)=(v(t),c(t))$.
\end{lemma}
\proof
Let $t \in (-1,1)$.
In the case of $\norm{\bm{v}(t)}_{E_\kappa}\geq 1$, the inequality \eqref{lin-con-2} follows Lemma \ref{lem-lin-est} under the assumption \eqref{ass-small-2}.
We consider the case of $ \norm{\bm{v}(t)}_{E_\kappa}<1$.
By Theorem \ref{thm-w-gwp} and the assumption \eqref{ass-lin-2}, we have
\begin{equation}\label{eq-lip-est-1}
\norm{P_+v(t)}_{E_{\kappa}}\leq \max\{e^{k_*t}, e^{k^*t}\}l_1 \norm{\bm{v}(0)}_{E_{\kappa}}^{\alpha}+C_\kappa \min\{\norm{\bm{v}(0)}_{E_{\kappa}}^2,\delta^2\}
\end{equation}
and
\begin{align}\label{eq-lip-est-2}
\norm{\bm{v}(t)}_{E_\kappa}^2 \geq & \norm{(P_0+P_-)v(0)}_{E_\kappa}^2+\min\{e^{2k^*t},e^{2k_*t}\}\norm{P_+v(0)}_{E_\kappa}^2+\norm{P_{\gamma}v(0)}_{E_\kappa}^2 \notag \\
&+|\log c(0)-\log c^*|^2-C\kappa \norm{v(0)}_{E_{\kappa}}^2 -C_\kappa \min\{ \norm{\bm{v}(0)}_{E_\kappa}^3,\delta^3\}.
\end{align}
% We define $l_*$ by 
% \begin{align*}
% \norm{P_+v(0)}_{E_\kappa}= l_* \norm{\bm{v}(0)}_{E_{\kappa}}^\alpha
% \end{align*}
By the inequality \eqref{eq-lip-est-2} and the concavity of $x^{\alpha/2}$, we obtain
\begin{equation}\label{eq-lip-est-3}
\norm{\bm{v}(t)}_{E_\kappa}^{\alpha} \geq \norm{\bm{v}(0)}_{E_\kappa}^{\alpha} -C(l_1^{\alpha}\norm{\bm{v}(0)}_{E_\kappa}^{\alpha} +(c\kappa+ C_\kappa \delta) \norm{\bm{v}(0)}_{E_\kappa}^{\alpha})
\end{equation}
for $|t|<1$.
Combining \eqref{eq-lip-est-1} and \eqref{eq-lip-est-3}, we have there exists $T^*>0$ such that
\begin{align*}
\norm{\bm{v}(t)}_{E_\kappa}^{\alpha} \geq& l_1^{-1}(\max\{e^{k_*t}, e^{k^*t}\}+l_1^{-1}C_\kappa \delta^{2-\alpha})^{-1}(1-Cl_1^{\alpha}-C\kappa^{1/2}-C_\kappa \delta^{1/2})\norm{P_+v(t)}_{E_{\kappa}}\\
\geq &
\begin{cases}
Cl_1^{-1} \norm{P_+v(t)}_{E_{\kappa}}, &|t|\leq T^*, \\
l_1^{-1} \norm{P_+v(t)}_{E_{\kappa}},& -T^*\leq t \leq -\frac{T^*}{2},
\end{cases}
\end{align*}

\qed

The following lemma shows that the flow map of the system \eqref{Leq-1}--\eqref{Leq-2} preserves the Lipschitz coefficients and the H\"older coefficients of graphs.

\begin{lemma}\label{lem-ind}
Let $1<\alpha<2$.
Assume that $l_1,l_2,\kappa,\delta>0$ satisfy the conditions $\eqref{ass-small}$ and $\eqref{ass-small-2}$.
There exists $T^*>0$ such that the solution map $U_{\delta}(t)$ of the system $\eqref{Leq-1}$--$\eqref{Leq-2}$ for $|t|\leq T^*$ defines a map $\mathcal{U}_{\delta}(t):\mathscr{G}^+_{l_1,l_2,\alpha,\delta,\kappa} \to \mathscr{G}^+_{C_\alpha l_1,C_Ll_2,\alpha,\delta,\kappa}$ uniquely by the relation $U_{\delta}(t)(\gbr{G}\times \R)=\gbr{\mathcal{U}_{\delta}(t)G}\times \R$.
Moreover, if $-T^*\leq t \leq -T^*/2$, then $\mathcal{U}_{\delta}(t)$ maps $\mathscr{G}^+_{l_1,l_2,\alpha,\delta,\kappa}$ into itself.
\end{lemma}
Lemma \ref{lem-ind} follows Lemma \ref{lem-lin-est} and the proof of Lemma 3.4 in \cite{N S 2}.

We define
\[\norm{G}_{\mathscr{G}^+_\kappa}=\sup_{\bm{v}\in \mathcal{H} \setminus \{(0,c^*)\}} \frac{\norm{G(\bm{v})}_{E_\kappa}}{\norm{\bm{v}}_{E_\kappa}}.\]
Since
\[\norm{G(\bm{v})}_{E_\kappa}\leq  l_2\norm{\bm{v}}_{E_\kappa}\]
for $G \in \mathscr{G}^+_{l_1,l_2,\alpha,\delta,\kappa}$ and $\bm{v}\in H^1(\RTL)\times (0,\infty)$,
we have $\norm{G}_{\mathscr{G}^+_\kappa}\leq l_2$ for $G \in \mathscr{G}^+_{l_1,l_2,\alpha,\delta,\kappa}$.
Therefore, $(\mathscr{G}^+_{l_1,l_2,\alpha,\delta,\kappa},\norm{\cdot}_{\mathscr{G}^+_{\kappa}})$ is a bounded complete metric space.

\begin{lemma}\label{lem-cont}
Under the conditions $\eqref{ass-small}$ and $\eqref{ass-small-2}$, the mapping $\mathcal{U}_{\delta}(t)$ is a contraction on $(\mathscr{G}^+_{l_1,l_2,\alpha,\delta,\kappa},\norm{\cdot}_{\mathscr{G}^+_{\kappa}})$ for $t<-T^*/2$.
\end{lemma}

% \proof

In the case $c^* \notin \{\frac{4n^2}{5L^2}; n \in \Z, n>1\}$, the proof of Lemma \ref{lem-cont} follows the proof of Lemma 4.3 in \cite{YY5}.
In the case $c^* \in \{\frac{4n^2}{5L^2};n\in \Z, n>1\}$, for $|T|<1$ there exists $C>0$ such that 
\[\norm{P_ae^{-T\mathcal{A}}v}_{E_\kappa} \leq \norm{P_av}_{E_\kappa} + C\kappa \norm{(I-P_1-P_a)v}_{E_{\kappa}} \leq (1+C\kappa) \norm{v}_{E_{\kappa}}.\]
Therefore, the proof of Lemma 4.3 in \cite{YY5} yields the conclusion of Lemma \ref{lem-cont} for $c^* \in \{\frac{4n^2}{5L^2};n\in \Z, n>1\}$ and sufficiently small $\kappa>0$.

From Lemma \ref{lem-cont}, we obtain the fix point of $\mathcal{U}_{\delta}$.
\begin{proposition}\label{prop-fixed-p}
Let $1<\alpha<2$. Assume that $l_1,l_2,\delta,\kappa>0$ satisfy the conditions $\eqref{ass-small}$ and $\eqref{ass-small-2}$.
There exists a unique $G_{+}^{\delta} \in \mathscr{G}^+_{l_1,l_2,\alpha,\delta,\kappa}$ such that $\mathcal{U}_{\delta}(t)G_+^{\delta}=G_+^{\delta}$ for all $t<0$, where $G_+^{\delta}$ does not depend on $\kappa$.
Moreover, the uniqueness hold for any fixed $t<0$.
\end{proposition}
\proof
The proof of the existence and the uniqueness for $t<0$ follows the proof of Theorem 3.6 in \cite{N S 2}.

We show that $G_+^{\delta}$ does not depend on $\kappa$.
Let $G_+^{\delta,\kappa}$ be the fixed point of $\mathcal{U}_{\delta}$ on $\mathscr{G}^+_{l_1,l_2,\alpha,\delta,\kappa}$.
By the definition of $\norm{\cdot}_{E_\kappa}$ we have $\mathscr{G}^+_{l_1,l_2,\alpha,\delta,\kappa_1} \subset \mathscr{G}^+_{l_1,l_2,\alpha,\delta,\kappa_2}$ for $\kappa_1 <\kappa_2$.
% \[ \norm{G(\bm{v}_0)-G(\bm{v}_1)}_{E_{\kappa_1}}=\norm{G(\bm{v}_0)-G(\bm{v}_1)}_{E_{\kappa_2}} \leq l \mathfrak{m}_{\delta,\kappa_2}(\bm{v}_0,\bm{v}_1) \leq l\mathfrak{m}_{\delta,\kappa_1}(\bm{v}_0,\bm{v}_1)\]
% for $G:\mathcal{H} \to P_+H^1(\RTL)$, $\kappa_1>\kappa_2$ and $\bm{v}_0,\bm{v}_1 \in \mathcal{H}$.
Therefore, by the uniqueness of the fixed point $G_+^{\delta,\kappa_2}$ of $\mathcal{U}_{\delta}$ on $\mathscr{G}^+_{l_1,l_2,\delta,\kappa_2}$, we have $G_+^{\delta,\kappa_2}=G_+^{\delta,\kappa_1}$ for $(l_1,l_2,\alpha,\delta,\kappa_1)$ and $(l_1,l_2,\alpha,\delta,\kappa_2)$ satisfying $\eqref{ass-small}$ and $\kappa_1<\kappa_2$.

\qed

In the following proposition, we show the existence of solitary waves near by $Q_{c^*}$ for $c^* \in \{\frac{4n^2}{5L^2};n\in \Z, n>1\}$. 
\begin{proposition}\label{prop-bifur}
Let $c^* =\frac{4n_0^2}{5L^2}$ for some $n_0 \in \Z \setminus\{-1,0,1\}$.
There exist $C_*,\delta_0>0$, $\check{c} \in C^2((-\delta_0,\delta_0)^2,(0,\infty))$ and for $\varphi_{c^*} \in C^2((-\delta_0,\delta_0)^2,H^2(\RTL))$ such that $\varphi_{c^*}(\bm{a})>0$, $\varphi_{c^*}(\bm{a})(x,y)=\varphi_{c^*}(\bm{a})(-x,y)$, $\check{c}(\bm{a})=\check{c}(|\bm{a}|,0)$,
\[-\Delta \varphi_{c^*}(\bm{a})+\check{c}(\bm{a})\varphi_{c^*}(\bm{a})-(\varphi_{c^*}(\bm{a}))^2=0,\]
\[\varphi_{c^*}(\bm{a})=Q_{c^*}+a_0Q_{c^*}^{3/2}\cos \frac{n_0y}{L} +a_1 Q_{c^*}^{3/2}\sin \frac{n_0y}{L} +O(|\bm{a}|^2) \mbox{ as } |\bm{a}| \to \infty,\]
\[\norm{\varphi_{c^*}(\bm{a})}_{L^2}^2=\norm{Q_{c^*}}_{L^2}^2+\frac{C_{2,c^*}}{2}|\bm{a}|^2+o(|\bm{a}|^2) \mbox{ as } |\bm{a}| \to \infty,\]
and
\[\check{c}(\bm{a})=c^*+\frac{C_*}{2}|\bm{a}|^2+o(|\bm{a}|^2)  \mbox{ as } |\bm{a}| \to \infty\]
for $\bm{a}=(a_0,a_1) \in (-\delta_0,\delta_0)^2$, where 
\[C_{2,c^*}=\frac{3C_*\norm{Q_{c^*}}_{L^2}^2}{2c^*}-\frac{5\norm{Q_{c^*}^{3/2}\cos \frac{n_0y}{L}}_{L^2}^2}{2}>0.\] 
\end{proposition}
In the case $n_0=1$, Proposition \ref{prop-bifur} was shown in Proposition 1.3 of \cite{YY4}.
By applying the Lyapunov--Schmidt reduction and the Crandall--Rabinowitz Transversality in \cite{C R,K K P} and the proof of Proposition 1.3 in \cite{YY4}, we can prove Proposition \ref{prop-bifur} for $n_0 \in \Z \setminus\{-1,0,1\}$.

We define $\mathfrak{g}(w,c)$ by
\[\mathfrak{g}(w,c)=w+G_+^{\delta}(w,c)+Q_c\]
for $w \in H^1(\RTL), c >0$.
Let 
\begin{align*}
\mathcal{M}_{cs}^{\delta}(c^*,r)=\{\tau_{\rho}\mathfrak{g}(w,c)&;w \in (P_-+P_a+P_{\gamma})H^1(\RTL),|c-c^*|<c^*/2, \\
& \inf_{q\in \R} \norm{\mathfrak{g}(w,c)-\tau_qQ_{c^*}}_{H^1}<r,\rho \in \R\}
\end{align*}
and
\begin{align*}
\tilde{\mathcal{M}}_{cs}^{\delta}(c^*,r)=\{\tau_{\rho}\mathfrak{g}(w,c)&; w \in P_{\leq 0}H^1(\RTL), \norm{(P_1+P_2)w}_{H^1(\RTL)}<r^{1/2},\\
& |c-c^*|<c^*/2,\inf_{q\in \R}\norm{\mathfrak{g}(w,c)-\tau_qQ_{c^*}}_{H^1}<r,\rho \in \R\}
\end{align*}
for $r>0$.
%By definition of $G^+_{\delta}$, $\mathcal{M}_{cs}^{\delta}(c^*,r)$ and $\tilde{\mathcal{M}}_{cs}^{\delta}(c^*,r)$ is invariant with respect to the flow of \eqref{LZKeq-1}--\eqref{LZKeq-2}}.

We show the stability of $Q_{c^*}$ on $\tilde{\mathcal{M}}_{cs}^{\delta}$ by using the conservation law $S_c(u)=E(u)+\frac{c}{2}M(u)$ as a Lyapunov function.
To recover the degeneracy of the Lyapunov function $S_c(u)$ around $Q_{c^*}$ for $c^* \in \{\frac{4n^2}{5L^2};n\in \Z, n>1\}$,
we define the modulated solitary wave 
\[\Theta(\bm{a},c)=\frac{c}{c^*}\varphi_{c^*}(\bm{a})\Bigl( \sqrt{\frac{c}{c^*}}x,y\Bigr)\]
and the modulation 
\[\beta(\bm{a},c)=\frac{c^*\norm{Q_c}_{L^2}^{4/3}}{\norm{\varphi_{c^*}(\bm{a})}_{L^2}^{4/3}}\]
for $\bm{a}\in \R^2$ and $c>0$.
Then, 
\begin{equation}\label{eq-beta-1}
\norm{\Theta(\bm{a},\beta(\bm{a},c))}_{L^2}=\norm{Q_{c}}_{L^2}
\end{equation}
 and
\begin{equation}\label{eq-beta-2}
\beta(\bm{a},c)-c=-c\frac{\norm{\varphi_{c^*}(\bm{a})}_{L^2}^{4/3}-\norm{Q_{c^*}}_{L^2}^{4/3}}{\norm{\varphi_{c^*}(\bm{a})}_{L^2}^{4/3}}
=-\frac{cC_{2,c^*}|\bm{a}|^2}{3\norm{Q_{c^*}}_{L^2}^2}+o(|\bm{a}|^2).
\end{equation}

In the following lemma, we investigate the fourth order term of the Lyapunov function.
\begin{lemma}\label{lem-Lyap}
For $c^* \in \{\frac{4n^2}{5L^2};n\in \Z,n>1\}$, $c>0$ and $\bm{a} \in \R^2$,
\begin{align*}
S_c(\Theta(\bm{a},\beta(\bm{a},c)))-S_c(Q_c)=&\Bigl(\frac{c}{c^*} \Bigr)^{5/2}\frac{5c^*C_{2,c^*}\norm{Q_{c^*}^{3/2}\cos \frac{y}{L}}_{L^2}^2|\bm{a}|^4}{48\norm{Q_{c^*}}_{L^2}^2}\\
&+\frac{c^*-c}{c^*}\norm{\partial_y\Theta(\bm{a},\beta(\bm{a},c))}_{L^2}^2+o(|\bm{a}|^4).
\end{align*}
\end{lemma}
Since $\norm{Q_{c^*}^{3/2}\cos \frac{ny}{L}}_{L^2}=\norm{Q_{c^*}^{3/2}\cos \frac{y}{L}}_{L^2} $ for non-zero integer $n$, the proof of Lemma \ref{lem-Lyap} is same as the proof of Lemma 5.3 in \cite{YY4}.
Therefore, we omit the proof of Lemma \ref{lem-Lyap}.

We define the orthogonality condition 
\begin{align}\label{eq-cri-mod-1}
(v,\Theta(\bm{a},c))_{L^2}=(v,\partial_x\Theta(\bm{a},c))_{L^2}=&(v,\partial_{a_0}\Theta(\bm{a},c))_{L^2}\notag\\
=&(v,\partial_{a_1}\Theta(\bm{a},c))_{L^2}=0.
\end{align}
for $(v,c,\bm{a})\in L^2(\RTL)\times (0,\infty) \times \R^2$
\begin{lemma}\label{lem-cri-mod}
Let $c^* \in \{\frac{4n^2}{5L^2};n\in \Z,n>1\}$.
 There exist $\delta, C>0$, $c:\mathcal{N}_{\delta,c^*} \to (0,\infty)$, $\rho:\mathcal{N}_{\delta,c^*}\to \R$ and $\bm{a}=(a_0,a_1):\mathcal{N}_{\delta,c^*} \to \R^2$ such that for $u \in \mathcal{N}_{\delta,c^*}$, $(v(u),c(u),\bm{a}(u))$ satisfies the orthogonality condition $\eqref{eq-cri-mod-1}$ and
\begin{equation}\label{eq-cri-mod-0}
\norm{v(u)}_{H^1}+|c(u)-c^*|+|\bm{a}(u)| \leq C \inf_{q \in \R}\norm{u-\tau_qQ_{c^*}}_{H^1}, 
\end{equation}
where
 \[v(u)=\tau_{-\rho(u)}u-\Theta(\bm{a}(u),c(u)).\]
\end{lemma}
\proof
We define
\[F_q(u,c,\rho,\bm{a})=
\begin{pmatrix}
(\tau_{-\rho}u-\tau_q\Theta(\bm{a},c),\tau_q\Theta(\bm{a},c))_{L^2}\\
(\tau_{-\rho}u-\tau_q\Theta(\bm{a},c),\tau_q\partial_x\Theta(\bm{a},c))_{L^2}\\
(\tau_{-\rho}u-\tau_q\Theta(\bm{a},c),\tau_q\partial_{a_1}\Theta(\bm{a},c))_{L^2}\\
(\tau_{-\rho}u-\tau_q\Theta(\bm{a},c),\tau_q\partial_{a_2}\Theta(\bm{a},c))_{L^2}
\end{pmatrix}\]
for $c>0, q, \rho \in \R, \bm{a} \in \R^2$ and $u \in H^1(\RTL)$.
Then, $F_q(\tau_qQ_{c^*},c^*,0,(0,0))= {}^t(0,0,0,0)$ and
\begin{align*}
&\frac{\partial F_q}{\partial c \partial \rho \partial a_1 \partial a_2}|_{(u,c,\rho,\bm{a})=(\tau_qQ_{c^*},c^*,0,(0,0))}\\
=&\mbox{diag}(-(Q_{c^*},\partial_cQ_{c^*})_{L^2},-\norm{\partial_xQ_{c^*}}_{L^2}^2,-\norm{Q_{c^*}^{3/2}\cos \frac{n_0y}{L}}_{L^2}^2, -\norm{Q_{c^*}^{3/2}\sin \frac{n_0y}{L}}_{L^2}^2).
\end{align*}
By the implicit function theorem, there exists $\delta, C>0$, $c_q:\mathcal{B}_{\delta}(\tau_qQ_{c^*}) \to (0,\infty)$, $\rho_q:\mathcal{B}_{\delta}(\tau_qQ_{c^*})\to \R$ and $\bm{a}_q:\mathcal{B}_{\delta}(\tau_qQ_{c^*}) \to \R^2$ such that $(\tau_{-q}v_q(u),c_q(u),\bm{a}_q(u))$ satisfy the orthogonality condition \eqref{eq-cri-mod-1},
\[\norm{v_q(u)}_{H^1}+|c_q(u)-c^*|+|\bm{a}_q(u)| \leq C \norm{u-\tau_qQ_{c^*}}_{H^1}\]
and
\begin{equation}\label{eq-cri-mod-2}
(c_{q_1}(u),\rho_{q_1}(u),\bm{a}_{q_1}(u))=(c_{q_2}(u),q_1-q_2+\rho_{q_2}(u),\bm{a}_{q_2}(u)), \quad u \in \mathcal{B}_{\delta}(\tau_{q_1}Q_{c^*})\cap\mathcal{B}_{\delta}(\tau_{q_2}Q_{c^*}),
\end{equation}
where $v_q(u)=\tau_{-\rho(u)}u-\tau_q\Theta(\bm{a}(u),c(u))$ and $\mathcal{B}_{\delta}(f)=\{u \in H^1; \norm{u-f}_{H^1}<\delta\}$.
By the compatibility condition \eqref{eq-cri-mod-2}, we can define $(c,\rho,\bm{a}):\mathcal{N}_{\delta,c^*} \to (0,\infty) \times \R \times \R^2$ by
\[(c(u),\rho(u),\bm{a}(u))=(c_q(u),\rho_q(u)-q,\bm{a}_q(u)), \quad u \in \mathcal{B}_{\delta}(\tau_{q}Q_{c^*}).\]
Then, $(v(u),c(u),\bm{a}(u))$ satisfies the orthogonality condition $\eqref{eq-cri-mod-1}$ and \eqref{eq-cri-mod-0}.
\qed

In the following lemma, we show the estimate of the difference the modulated solitary waves.
\begin{lemma}\label{lem-est-theta}
Let $c^* \in \{\frac{4n^2}{5L^2};n\in \Z,n>1\}$.
There exists $\delta>0$ such that for $c_0>0$ and $u\in \mathcal{N}_{\delta,c^*}$ with $\norm{u}_{L^2}=\norm{Q_{c_0}}_{L^2}$ and $|c_0-c^*|\ll \delta^{1/2}$,
\[ \norm{\Theta(\bm{a}(u), \beta(\bm{a}(u),c_0))-\Theta(\bm{a}(u),c(u))}_{H^1} \lesssim \norm{v(u)}_{L^2}^2,\]
\[|\beta(\bm{a}(u),c_0)-c(u)|\lesssim \norm{v(u)}_{L^2}^2,\]
where $\bm{a}(u),c(u)$ and $v(u)$ is defined in Lemma $\ref{lem-cri-mod}$.
\end{lemma}
\proof
For $u \in \mathcal{N}_{\delta,c^*}$ with $\norm{u}_{L^2}=\norm{Q_{c_0}}_{L^2}$, the equation \eqref{eq-beta-1} yields 
\begin{align*}
\norm{\Theta(\bm{a}(u),\beta(\bm{a}(u),c_0))}_{L^2}^2=\norm{Q_{c_0}}_{L^2}^2=&\norm{v(u)+\Theta(\bm{a}(u),c(u))}_{L^2}^2\\
=&\norm{v(u)}_{L^2}^2+\norm{\Theta(\bm{a}(u),c(u))}_{L^2}^2.
\end{align*}
Since
\[|c(u)-c^*|+|\beta(\bm{a}(u),c_0)-c^*|<\frac{c^*}{2}\]
for sufficiently small $\delta>0$, we have
\begin{align*}
\norm{v(u)}_{L^2}^2=&\norm{\Theta(\bm{a}(u),\beta(\bm{a}(u),c_0))}_{L^2}^2-\norm{\Theta(\bm{a}(u),c(u))}_{L^2}^2\\
=&(\beta(\bm{a}(u),c_0)^{\frac{3}{2}}-c(u)^{\frac{3}{2}})\norm{\varphi_{c^*}(\bm{a}(u))}_{L^2}^2
\gtrsim  \beta(\bm{a}(u),c_0)-c(u)\geq 0.
\end{align*}
Therefore, we obtain 
\begin{align*}
&\norm{\Theta(\bm{a}(u),\beta(\bm{a}(u),c_0))-\Theta(\bm{a}(u),c(u))}_{H^1}\\
\lesssim & |\beta(\bm{a}(u),c_0)-c(u)|\norm{\partial_cQ_{c^*}}_{H^1}+o(|\beta(\bm{a}(u),c_0)-c(u)|)
\lesssim  \norm{v(u)}_{L^2}^2.
\end{align*}

\qed

In the following theorem, we prove the stability of the line solitary wave on the set $\tilde{\mathcal{M}}_{cs}^{\delta}(c^*,\varepsilon)$.
\begin{theorem}\label{thm-stability}
Let $3/2<\alpha<2$.
Assume that $l_1,l_2,\delta,\kappa>0$ satisfy the condition $\eqref{ass-small}$ and $\eqref{ass-small-2}$.
For any $\varepsilon >0$, there exists $\tilde{\varepsilon}=\tilde{\varepsilon}(c^*,\varepsilon)>0$ such that for $u_0 \in \tilde{\mathcal{M}}_{cs}^{\delta}(c^*,\tilde{\varepsilon})$ the solution $u$ to the equation $\eqref{ZKeq}$ with the initial data $u_0$ satisfies $u(t) \in \tilde{\mathcal{M}}_{cs}^{\delta}(c^*,\varepsilon)$ for all $t>0$.
\end{theorem}
\proof
In the case $c^* \notin \{\frac{4n^2}{5L^2};n\in \Z, n>1\}$, the conclusion follows the proof of Theorem 4.6 in \cite{YY5}.

We show the case $c^* \in \{\frac{4n^2}{5L^2};n\in \Z, n>1\}$.
Let $l_1,l_2,\delta,\kappa>0$ satisfying \eqref{ass-small} and \eqref{ass-small-2}.
We prove the conclusion by contradiction.
We assume there exists $0<\varepsilon_0 \ll \delta^2$ such that for $0<\tilde{\varepsilon}<\varepsilon_0$ there exist $t_0>0$ and a solution $u$ to $\eqref{ZKeq}$ with an initial data $\tau_{\rho_0}\mathfrak{g}(w_0,c_0) \in \tilde{\mathcal{M}}_{cs}^{\delta}(c^*,\tilde{\varepsilon})$ satisfying $\norm{(P_1+P_2)w_0}_{H^1}\leq \tilde{\varepsilon}$,  
\[
\sup_{0\leq t \leq t_0}\inf_{q \in \R}\norm{u(t)-\tau_qQ_{c^*}}_{H^1}\leq \varepsilon_0
\]
and
\begin{equation}\label{eq-stabiltiy-ass}
\inf_{q\in \R}\norm{u(t_0)-\tau_qQ_{c^*}}_{H^1}=\varepsilon_0.
\end{equation}
We define the solution $(v_1(t),c_1(t),\rho_1(t))$ to the system \eqref{Leq-1}--\eqref{Leq-2} with the initial data $(w_0+G_+^{\delta}(w_0,c_0),c_0,\rho_0)$ and $c_2(t)=c(u(t))$, $\rho_2(t)=\rho(u(t))$, $\bm{a}_2(t)=\bm{a}(u(t))$ and 
\[v_2(t)=\tau_{\rho_1(t)-\rho_2(t)}(v_1(t)+Q_{c_1(t)})-\Theta(\bm{a}_2(t),c_2(t)),\]
where $c(u), \rho(u)$ and $\bm{a}(u)$ are defined in Lemma \ref{lem-cri-mod}.
Then, $(v_2(t),c_2(t),\bm{a}_2(t))$ satisfies the orthogonality condition \eqref{eq-cri-mod-1}. 
Since
\begin{align}%\label{eq-stabitliy-1}
\norm{\tau_{\rho_1}Q_{c}-\tau_qQ_{c^*}}_{H^1}\lesssim & \norm{(P_1+P_2)(Q_{c}-\tau_{q-\rho_1}Q_{c^*})}_{H^1}\notag\\
\lesssim & \norm{(P_1+P_2)v}_{H^1}+\norm{\tau_{\rho_1}(v+Q_{c})-\tau_qQ_{c^*}}_{H^1}, \notag
\end{align}
we obtain
\begin{equation*}
\norm{v_1(t)}_{H^1}+|c_1(t)-c^*|\lesssim \inf_{q \in \R} \norm{\tau_{\rho_1(t)}(v_1(t)+Q_{c_1(t)})-\tau_{q}Q_{c^*}}_{H^1}+\norm{(P_1+P_2)v_1(0)}_{H^1}.
\end{equation*}
Therefore, by the continuity of $u(t)$ and $(v_1(t),c_1(t),\rho_1(t))$ we have
\[\norm{v_1(t)}_{H^1}+|c_1(t)-c^*|\lesssim \varepsilon_0+\tilde{\varepsilon}^{1/2} \ll \delta\]
and
\begin{equation}\label{eq-stability-2-1}
u(t)=\tau_{\rho_1(t)}(v_1(t)+Q_{c_1(t)})=\tau_{\rho_2(t)}\bigl(v_2(t)+\Theta(\bm{a}_2(t),c_2(t))\bigr)
\end{equation}
for $0 \leq t \leq t_0$.
%By definition of $G_+^{\delta}$, we have $P_+(v_1(t))=G_+^{\delta}(v_1(t),c_1(t))$ and $u(t)=\tau_{\rho_1(t)}\mathfrak{g}(P_{\leq 0} v_1(t),c_1(t)) \in \tilde{\mathcal{M}}_{cs}^{\delta}(c^*,\varepsilon_0)$ for $0 \leq t \leq t_0$.
We define $c_u>0$ by $\norm{u(0)}_{L^2}=\norm{Q_{c_u}}_{L^2}$.
By Lemma \ref{lem-est-theta}, we have
\begin{equation}\label{eq-stability-2}
\norm{\Theta(\bm{a}_2(t),\beta(\bm{a}_2(t),c^*))-\Theta(\bm{a}_2(t),c_2(t))}_{H^1}\lesssim |c^*-c_u|+\norm{v_2(t)}_{L^2}^2.
\end{equation}
Therefore, we obtain
\begin{align}\label{eq-stability-3}
&\tbr{S_{c^*}'\bigl(\Theta(\bm{a}_2,\beta(\bm{a}_2,c^*))\bigr),v_2}_{H^{-1},H^1}\notag \\
=&\tbr{S_{\check{c}(\bm{a}_2)}'(\varphi_{c^*}(\bm{a}_2)),v_2}_{H^{-1},H^1}+\tbr{S_{\check{c}(\bm{a}_2)}''(\varphi_{c^*}(\bm{a}_2))(\Theta(\bm{a}_2,\beta(\bm{a}_2,c^*))-\varphi_{c^*}(\bm{a}_2)),v_2}_{H^{-1},H^1}\notag \\
&+(c^*-\check{c}(\bm{a}_2))\tbr{\Theta(\bm{a}_2,\beta(\bm{a}_2,c^*))-\Theta(\bm{a}_2,c_2),v_2}_{H^{-1},H^1}\notag\\
&+o(\norm{\Theta(\bm{a}_2,\beta(\bm{a}_2,c^*))-\varphi_{c^*}(\bm{a}_2)}_{H^1}\norm{v_2}_{H^1})\notag\\
=&o((|c_u-c^*|+|\beta(\bm{a}_2,c^*)-c^*|+\norm{v_2}_{L^2})\norm{v_2}_{L^2})
\end{align}
Thus, from \eqref{eq-beta-2}, \eqref{eq-cri-mod-1}, \eqref{eq-stability-2} and \eqref{eq-stability-3} there exists $C>0$ such that
\begin{align}\label{eq-stability-4}
&S_{c^*}(u)-S_{c^*}(Q_{c^*})\notag\\
=&S_{c^*}\bigl(\Theta(\bm{a}_2,\beta(\bm{a}_2,c^*))\bigr)-S_{c^*}(Q_{c^*})+\frac{1}{2}\tbr{S_{c^*}''(Q_{c^*})(I-P_0)v_2,(I-P_0)v_2}_{H^{-1},H^1} \notag \\
&+o(\norm{v_2}_{L^2}^2+|c^*-c_u|+|\bm{a}_2|^4)\notag \\
\geq& C|\bm{a}_2|^4 + \frac{1}{2}\tbr{\mathbb{L}_{c^*}\gamma_2,\gamma_2}_{H^{-1},H^1}+\sum_{j,k}\Lambda_{k,2}^{+,j}\Lambda_{k,2}^{-,j}+o(\norm{v_2}_{H^1}^2+|c^*-c_u|+|\bm{a}_2|^4)
\end{align}
where 
\[\Lambda_{k,2}^{\pm,j}=(v_2,\mathbb{L}_{c^*}F_k^{\mp,j})_{L^2}, \quad \gamma_2=P_{\gamma}v_2\]
for $j\in \{0,1\}$.
By Lemma \ref{lem-cri-mod}, we have
\begin{equation}\label{eq-stability-4-2}
\sup_{0\leq t \leq t_0}(\norm{v_2(t)}_{H^1}+|c_2(t)-c^*|+|\bm{a}_2(t)|)\lesssim \varepsilon_0.
\end{equation}
The equation \eqref{eq-stability-2-1} yields
\begin{align*}
&|(\tau_{\rho_1-\rho_2}Q_{c^*},\partial_xQ_{c^*})_{L^2}|\\
\leq & |(\tau_{\rho_1-\rho_2}Q_{c_1},\partial_x\Theta(\bm{a}_2,c_2))_{L^2}-(\tau_{\rho_1-\rho_2}Q_{c^*},\partial_xQ_{c^*})_{L^2}|+|(v_1,\partial_x Q_{c^*})_{L^2}| \\
&+|(v_1,\tau_{\rho_2-\rho_1}\partial_x\Theta(\bm{a}_2,c_2))_{L^2}-(v_1,\partial_x Q_{c^*})_{L^2}|\\
\lesssim & (\varepsilon_0+|\rho_1-\rho_2|)^2+\tilde{\varepsilon}^{1/2}
\end{align*}
for sufficiently small $|\rho_1-\rho_2|$.
Therefore, if $|\rho_1-\rho_2|$, $\varepsilon_0$ and $\tilde{\varepsilon}$ are sufficiently small, we have
\[|\rho_1-\rho_2|\lesssim |(\tau_{\rho_1-\rho_2}Q_{c^*},\partial_xQ_{c^*})_{L^2}| \lesssim \varepsilon_0^2+ \tilde{\varepsilon}^{1/2}\]
which yields 
\begin{equation*}
\sup_{0\leq t \leq t_0}|\rho_1(t)-\rho_2(t)|\lesssim \varepsilon_0^2 + \tilde{\varepsilon}^{1/2}.
\end{equation*}
Since $P_+v_1(t)=G_+^{\delta}(v_1(t),c_1(t))$, we obtain
\begin{align}\label{eq-stability-5}
\sum_{j,k}|\Lambda_{k,2}^{+,j}(t)|^2=&\norm{P_+v_2(t)}_{E_{\kappa}}^2\notag\\
\lesssim & \norm{P_+v_1(t)}_{E_{\kappa}}^2+|\rho_1(t)-\rho_2(t)|^2\norm{v_1(t)}_{E_{\kappa}}^2\notag\\
=&\norm{G_+^{\delta}(v_1(t),c_1(t))}_{E_{\kappa}}^2+|\rho_1(t)-\rho_2(t)|^2\norm{v_1(t)}_{E_{\kappa}}^2\notag\\
\lesssim & l_1^2(\norm{v_1(t)}_{E_{\kappa}}^2+|c_1(t)-c^*|^2)^{\alpha} +|\rho_1(t)-\rho_2(t)|^2\norm{v_1(t)}_{E_{\kappa}}^2 \notag \\
\lesssim &l_1^2(\varepsilon_0^2+\tilde{\varepsilon})^{\alpha}+(\varepsilon_0^4+\tilde{\varepsilon}^2).
\end{align}
% From the inequality \eqref{eq-stability-4} and \eqref{eq-stability-4-2}, we obtain 
% \begin{align}\label{eq-stability-4-1}
% |\bm{a}_2|^4+\norm{v_2}_{H^1}^2\lesssim \tilde{\varepsilon}^2+\varepsilon_0\sum_{j,k}|\Lambda_{k,2}^{-,j}|
% \end{align}
By the definition of $v_2$, we have
\begin{equation}\label{eq-stability-5-1}
\partial_t\Lambda_{k,2}^{-,j}=-\lambda_k\Lambda_{k,2}^{-,j}+\bigl(-\partial_x[(v_2)^2+2v_2(\Theta(\bm{a}_2,c_2)-Q_{c^*})],\mathbb{L}_{c^*}F_{k}^{+,j})_{L^2}.
\end{equation}
Since $|\Lambda_{k,2}^{-,j}(0)|\lesssim \tilde{\varepsilon}$, for sufficiently small $\tilde{\varepsilon}$ we have 
\begin{equation}\label{eq-stability-6}
\sup_{0\leq t \leq t_0} |\Lambda_{k,2}^{-,j}(t)| \lesssim \tilde{\varepsilon}+\varepsilon_0^2.
\end{equation}
Combining \eqref{eq-stability-4}--\eqref{eq-stability-6}, (v) of Lemma \ref{prop-linear} and Lemma \ref{lem-Lyap}, we obtain 
\begin{align}\label{eq-stability-7}
&\sup_{0\leq t \leq t_0} (\norm{(I-P_--P_+)v_2(t)}_{H^1}^2+|\bm{a}_2(t)|^4)\notag \\
\lesssim & S_{c^*}(u)-S_{c^*}(Q_{c^*})+l_1((\varepsilon_0+\tilde{\varepsilon}^{1/2})^{\alpha}+(\varepsilon_0^2+\tilde{\varepsilon}))(\varepsilon_0^2+\tilde{\varepsilon})\notag\\
\lesssim & l_1(\varepsilon_0^{\alpha+2}+\tilde{\varepsilon}).
\end{align}
Using the estimates \eqref{eq-stability-5}, \eqref{eq-stability-6} and \eqref{eq-stability-7}, we obtain
\begin{equation}\label{eq-stability-7-1}
\sup_{0\leq t \leq t_0} |\Lambda_{k,2}^{-,j}(t)| \lesssim \tilde{\varepsilon}+\varepsilon_0^{\alpha+1}.
\end{equation}
By the bootstrap argument, we have
\begin{equation}\label{eq-stability-7-2}
\sup_{0\leq t \leq t_0} (\norm{(I-P_--P_+)v_2(t)}_{H^1}^2+|\bm{a}_2(t)|^4)
\lesssim  l_1(\varepsilon_0^{2\alpha+1}+\tilde{\varepsilon}).
\end{equation}
% By the equation \eqref{eq-beta-2} and Lemma \ref{lem-est-theta}, we have
% \begin{align}\label{eq-stability-8}
% &\sup_{0\leq t \leq t_0} |c_2(t)-c^*| \leq \sup_{0\leq t \leq t_0}(|c_2(t)-\beta(\bm{a}_2(t),c_u)|+|\beta(\bm{a}_2(t),c_u)-c_u|+|c_u-c^*|)\notag \\
%  \lesssim& (l+\varepsilon_0+\tilde{\varepsilon}^{1/2})(\varepsilon_0+\tilde{\varepsilon}^{1/2}).
% \end{align}
For $\tilde{\varepsilon} \ll \varepsilon_0$, the inequalities \eqref{eq-stability-5}, \eqref{eq-stability-7-1} and \eqref{eq-stability-7-2} yield 
\begin{align*}
&\sup_{0\leq t \leq t_0} \inf_{q \in \R} \norm{u(t)-\tau_q Q_{c^*}}_{H^1} \\
\lesssim & \sup_{0\leq t \leq t_0} (\norm{v_2(t)}_{H^1}+\norm{\Theta(\bm{a}_2(t),c_2(t))-\Theta(\bm{a}_2(t),\beta(\bm{a}_2(t),c^*))}_{H^1}\\
&+\norm{\Theta(\bm{a}_2(t),\beta(\bm{a}_2(t),c^*))-\varphi_{c^*}(\bm{a}_2(t))}_{H^1}+\norm{\varphi_{c^*}(\bm{a}_2(t))-Q_{c^*}}_{H^1})\\
\lesssim & \varepsilon_0^{\alpha}+\varepsilon_0^{\frac{2\alpha+1}{4}}+\tilde{\varepsilon}^{1/2} \ll \varepsilon_0
\end{align*}
which contradicts the inequality \eqref{eq-stabiltiy-ass}.
Therefore, for any $\varepsilon>0$ there exists $\tilde{\varepsilon}>0$ such that for $u_0 \in \tilde{\mathcal{M}}_{cs}^{\delta}(c^*,\tilde{\varepsilon})$ the solution $u$ to the equation $\eqref{ZKeq}$ with the initial data $u_0$ satisfies 
\[ \sup_{ t\geq 0} \inf_{q \in \R} \norm{u(t)-\tau_qQ_{c^*}}_{H^1}<\varepsilon.\]
By the same calculation to obtain the equality \eqref{eq-stability-2-1} and $P_+(v_1(t))=G_+^{\delta}(v_1(t),c_1(t))$, we have $u(t)\in \tilde{\mathcal{M}}_{cs}^{\delta}(c^*,\varepsilon)$ for all $t>0$.
\qed
% Thus, we show the case $c \in \{\frac{4n^2}{5L^2};n\in \Z\}$

In the following lemma, we show a lower estimate for the unstable mode.
\begin{lemma}\label{lem-exit-est}
There exists $C_{\kappa}>0$ such that if $\kappa,\delta,l_0>0$ satisfy
\begin{equation}\label{ass-exit}
(\kappa^{1/2}+C_{\kappa}\delta^{1/4})(1+l_0)\ll \min\{1,l_0\},
\end{equation}
then there exists $T^*>0$ such that for any solutions $(v_0,c_0,\rho_0)$ and $(v_1,c_1,\rho_1)$ to the system $\eqref{Leq-1}$--$\eqref{Leq-2}$ satisfying
\begin{equation}\label{eq-exit-est-1}
\mathfrak{m}_{\delta,\kappa}(\bm{v}_0(0),\bm{v}_1(0))^2-\norm{P_+(v_0(0)-v_1(0))}_{E_\kappa}^2\leq l_0^2\norm{P_+(v_0(0)-v_1(0))}_{E_\kappa}^2,
\end{equation}
one has
\begin{equation*}
\mathfrak{m}_{\delta,\kappa}(\bm{v}_0(t),\bm{v}_1(t))^2-\norm{P_+(v_0(t)-v_1(t))}_{E_\kappa}^2\leq
\begin{cases}
2l_0^2 \norm{P_+(v_0(t)-v_1(t))}_{E_\kappa}^2, & 0 \leq t < T^*/2, \\
l_0^2 \norm{P_+(v_0(t)-v_1(t))}_{E_\kappa}^2, &  T^*/2 \leq t \leq T^*,
\end{cases}
\end{equation*}
and
\begin{equation*}
\norm{P_+(v_0(t)-v_1(t))}_{E_\kappa} \geq 
\begin{cases}
\frac{1}{2}e^{k_*t/2}\norm{P_+(v_0(0)-v_1(0))}_{E_\kappa}, & 0\leq t < T^*/2,\\
e^{k_*t/2}\norm{P_+(v_0(0)-v_1(0))}_{E_\kappa}, & T^*/2\leq t \leq T^*.
\end{cases}
\end{equation*}
where $k_*$ is defined by $\eqref{def-k}$.
\end{lemma}
\proof
By the assumption \eqref{ass-exit}, we have
\begin{equation}\label{eq-exit-est-2}
\norm{P_+(v_0(0)-v_1(0))}_{E_{\kappa}} \leq \mathfrak{m}_{\delta,\kappa}(\bm{v}_0(0),\bm{v}_1(0))\leq (1+l_0)\norm{P_+(v_0(0)-v_1(0))}_{E_{\kappa}}.
\end{equation}
By the above inequality and Lemma \ref{lem-est-mobile}, we obtain 
\begin{equation}\label{eq-exit-est-3}
\norm{P_+(v_0(0)-v_1(0))}_{E_{\kappa}} \leq e^{-k_* t} \norm{P_+(v_0(t)-v_1(t))}_{E_{\kappa}} + C_{\kappa}\delta^{\frac{1}{4}}\mathfrak{m}_{\delta,\kappa}(\bm{v}_0(0),\bm{v}_1(0)).
\end{equation}
Lemma \ref{lem-est-mobile} and the inequalities \eqref{lin-eq-5}, \eqref{eq-exit-est-2} and \eqref{eq-exit-est-3} yield 
\begin{align*}
&\mathfrak{m}_{\delta,\kappa}(\bm{v}_0(t),\bm{v}_1(t))^2-\norm{P_+(v_0(t)-v_1(t))}_{E_{\kappa}}^2\\
\leq & e^{-2k_*t}\norm{P_-(v_0(0)-v_1(0))}_{E_{\kappa}}^2+\norm{P_0(v_0(0)-v_1(0))}_{E_{\kappa}}^2+I(v_0(0),v_1(0))\\
&+|\log c_0(0) -\log c_1(0) |^2 +(C\kappa+C_{\kappa}\delta^{\frac{1}{2}})\mathfrak{m}_{\delta,\kappa}(\bm{v}_0(0),\bm{v}_1(0))^2\\
\leq & (l_0^2 +(C\kappa+C_{\kappa}\delta^{\frac{1}{2}})(1+l_0)^2)(e^{2k_*t}-C_{\kappa}\delta^{\frac{1}{2}}(1+l_0)^2)^{-1}\norm{P_+(v_0(t)-v_1(t))}_{E_{\kappa}}^2
\end{align*}
for $0 \leq t \leq T^*$.
\qed

To prove (v) of Theorem \ref{thm-main}, we show the following lemma.

\begin{lemma}\label{lem-exit-pro}
Let $\delta, l_0, \kappa>0$.
Suppose the assumption $\eqref{ass-exit}$.
There exists $0<\varepsilon_*=\varepsilon_*(c^*,\delta,l_0,\kappa)<\delta$ such that for any $0<\varepsilon<\varepsilon_*$ and solutions $u_0(t)$ and $u_1(t)$ to the equation $\eqref{ZKeq}$ satisfying 
\begin{equation}\label{eq-exi-1}
\sup_{t \geq 0} \inf_{q \in \R} \norm{u_1(t)-\tau_qQ_{c^*}}_{H^1}<\varepsilon, \quad \inf_{q \in \R} \norm{u_0(0)-\tau_qQ_{c^*}}_{H^1}<\varepsilon
\end{equation}
and
\begin{equation}\label{eq-exi-2}
\mathfrak{m}_{\delta,\kappa}(\bm{v}_0(0),\bm{v}_1(0))^2-\norm{P_+(v_0(0)-v_1(0))}_{E_{\kappa}}^2<l_0^2\norm{P_+(v_0(0)-v_1(0))}_{E_{\kappa}}^2,
\end{equation}
one has
\begin{equation}\label{eq-exi-3}
\inf_{q\in \R}\norm{u_0(t_0)-\tau_qQ_{c^*}}_{H^1}\geq \varepsilon
\end{equation}
for some $t_0>0$, where $(v_0(0),c_0(0),\rho_0(0))$ and $(v_1(0),c_1(0),\rho_1(0))$ satisfy
\begin{equation}\label{eq-exi-4}
u_j(0)=\tau_{\rho_j(0)}(v_j(0)+Q_{c_j(0)}), \quad |(v_1(0),\partial_xQ_{c^*})_{L^2}|+|(v_1(0),Q_{c^*})_{L^2}|<\varepsilon^{1/2}
\end{equation}
for $j=0,1$.
\end{lemma}

\proof
Let $v_j$ be the solutions to the system \eqref{Leq-1}--\eqref{Leq-2} with the initial data $(v_j(0),c_j(0),\rho_j(0))$.
We prove the inequality \eqref{eq-exi-3} by the contradiction.
Assume for any $0<\varepsilon_*\ll \delta^2$ there exist $0<\varepsilon <\varepsilon_*$ and solutions $u_0(t)$ and $u_1(t)$ to the equation \eqref{ZKeq} satisfying \eqref{eq-exi-1}, \eqref{eq-exi-2}, \eqref{eq-exi-4} and
\begin{equation}\label{eq-exi-contradiction}
\sup_{t\geq 0}\inf_{q\in \R} \norm{u_0(t)-\tau_qQ_{c^*}}_{H^1}<\varepsilon.
\end{equation}
By the inequality
\[ \norm{v_1(t)}_{H^1}+|c_1(t)-c^*|\lesssim \inf_{q \in \R}\norm{u_1(t)-\tau_qQ_{c^*}}_{H^1}+|(v_1(0),\partial_xQ_{c^*})_{L^2}|+|(v_1(0),Q_{c^*})_{L^2}|\ll \delta
\]
and $u_1(t)=\tau_{\rho_1(t)}(v_1(t)+Q_{c_1(t)})$ as long as $\norm{v_1(t)}_{H^1}^2+|c_1(t)-c^*|^2 < \delta^2$, we obtain
$u_1(t)=\tau_{\rho_1(t)}(v_1(t)+Q_{c_1(t)})$ for all $t\geq 0$.
From the inequality \eqref{eq-exi-2}, applying Lemma \ref{lem-exit-est} repeatedly, we have
\begin{equation}\label{eq-exi-5}
0<\frac{1}{2}e^{\frac{k_*t}{2}}\norm{P_+(v_0(0)-v_1(0))}_{E_{\kappa}}\leq \norm{P_+(v_0(t)-v_1(t))}_{E_{\kappa}}
\end{equation}
and
\begin{equation}\label{eq-exi-6}
(\mathfrak{m}_{\delta,\kappa}(\bm{v}_0(t),\bm{v}_1(t))^2-\norm{P_+(v_0(t)-v_1(t))}_{E_{\kappa}}^2)^{\frac{1}{2}}<2l_0\norm{P_+(v_0(t)-v_1(t))}_{E_{\kappa}}
\end{equation}
for all $t>0$.
Since 
\begin{align*}
 \norm{P_+(v_0(t)-v_1(t))}_{E_{\kappa}}\lesssim & \norm{P_+v_0(t)}_{H^1}+\norm{P_+v_1(t)}_{H^1}\\
\lesssim & \inf_{q \in \R} \norm{v_0(t)+Q_{c_0(t)}-\tau_{q}Q_{c^*}}_{H^1}+\inf_{q \in \R}\norm{v_1(t)+Q_{c_1(t)}-\tau_qQ_{c^*}}_{H^1},
\end{align*}
by the assumption \eqref{eq-exi-1} and the inequality \eqref{eq-exi-5} we have
\begin{equation}\label{eq-exi-7}
\frac{1}{2}e^{\frac{k_*t}{2}}\norm{P_+(v_0(0)-v_1(0))}_{E_{\kappa}} \leq \norm{P_+(v_0(t)-v_1(t))}_{E_{\kappa}} \lesssim \inf_{q\in \R} \norm{u_0(t)-\tau_qQ_{c^*}}_{H^1}+\varepsilon
\end{equation}
for all $t\geq 0$ and $0 < \varepsilon \ll \delta$ as long as $\norm{v_0(t)}_{H^1}^2+|c_0(t)-c^*|^2<\delta^2$.
By the assumption \eqref{eq-exi-1} and the inequality \eqref{eq-exi-6}, we obtain 
\begin{align}\label{eq-exi-8}
\norm{v_0(t)}_{H^1}^2+|c_0(t)-c^*|^2\lesssim & \mathfrak{m}_{\delta,\kappa}(\bm{v}_0(t),(0,c^*))^2\notag\\
\lesssim & \mathfrak{m}_{\delta,\kappa}(\bm{v}_0(t),\bm{v}_1(t))^2+\inf_{q \in \R}\norm{(I-P_0)v_1(t)+Q_{c_1(t)}-\tau_qQ_{c^*}}_{H^1}^2\notag\\
&+ |(v_1(t),\partial_xQ_{c^*})_{L^2}|^2+|(v_1(t),Q_{c^*})_{L^2}|^2\notag\\
\lesssim & (1+l_0)^2\norm{P_+(v_0(t)-v_1(t))}_{E_{\kappa}}^2+\varepsilon.
\end{align}
The inequality \eqref{eq-exi-5}, \eqref{eq-exi-7} and \eqref{eq-exi-8} contradict the assumption \eqref{eq-exi-contradiction} for sufficiently small $\varepsilon_*>0$.
Therefore, the proof was completed.

\qed

In the next corollary, we show (v) of Theorem \ref{thm-main} which means that solutions with the initial data off $\mathcal{M}_{cs}^{\delta}$ go out of the neighborhood of the line solitary waves.
\begin{corollary}\label{cor-exit}
Let $\delta, l_0, \kappa>0$.
Suppose the assumption $\eqref{ass-exit}$.
There exists $\varepsilon^*=\varepsilon^*(c^*,\delta,l_0,\kappa)>0$ such that for $u(0)\in \mathcal{N}_{\varepsilon^*, c^*}\setminus \mathcal{M}_{cs}^{\delta}(c^*,\varepsilon^*)$, the solution $u$ of the equation $\eqref{ZKeq}$ corresponding to the initial data $u(0)$ satisfies
\[\inf_{q \in \R}\norm{u(t_0)-\tau_qQ_{c^*}}_{H^1}\geq \varepsilon^*\]
for some $t_0 \geq 0$.
\end{corollary}

\proof
By the Lipschitz continuity of $G_+^{\delta}$, there exists $C>1$ such that for $0<\varepsilon<\min\{\tilde{\varepsilon},\varepsilon_*\}$ and $u  \in \mathcal{N}_{\varepsilon,c^*}$, $\tau_{\rho(u)}( (P_-+P_a+P_{\gamma})v+G_+^{\delta}((P_-+P_a+P_{\gamma})v,c(u))+Q_{c(u)}) \in \mathcal{N}_{C\varepsilon,c^*}$, where  $c(u)$ and $\rho(u)$ is defined by Lemma \eqref{lem-orth} and $v=\tau_{-\rho(u)}u-Q_{c(u)}$.
Let $0<\varepsilon<\min\{C^{-1}\tilde{\varepsilon},\varepsilon_*\}$ and $u(0) \in \mathcal{N}_{\varepsilon} \setminus \mathcal{M}_{cs}^{\delta}(c^*,\varepsilon)$ and $u$ be the solution to the equation \eqref{ZKeq} corresponding to the initial data $u(0)$, where $\varepsilon_*$ is defined by Lemma \ref{lem-exit-pro} and $\tilde{\varepsilon}=\tilde{\varepsilon}(c^*,\varepsilon^*)$ is defined by Theorem \ref{thm-stability}. 
We define the solution $(v_1,c_1,\rho_1)$ to the system \eqref{Leq-1}--\eqref{Leq-2} corresponding to the initial data $((P_-+P_a+P_{\gamma})v(0)+G_+^{\delta}((P_-+P_a+P_{\gamma})v(0),c(u(0))),c(u(0)),\rho(u(0)))$, where $v(0)=\tau_{-\rho(u(0))}u(0)-Q_{c(u(0))}$.
Since $(P_-+P_a+P_{\gamma})v(0)+G_+^{\delta}((P_-+P_a+P_{\gamma})v(0),c(u(0)))+Q_{c(u(0))}\in \tilde{\mathcal{M}}_{cs}^{\delta}(c^*,\tilde{\varepsilon})$, by Theorem \ref{thm-stability} we obtain $u_1(t)=\tau_{\rho_1(t)}(v_1(t)+Q_{c_1(t)})$ is a solution to \eqref{ZKeq} and
\begin{equation}\label{eq-cor-exi-1}
\sup_{t\geq 0}\inf_{q\in \R} \norm{u(t)-\tau_qQ_c}_{H^1} < \varepsilon_*.
\end{equation}
By the definition of $(v_1,c_1,\rho_1)$, we have the solution $(v_0,c_0,\rho_0)$ to the system \eqref{Leq-1}--\eqref{Leq-2} corresponding to the initial data $(v(0),c(u(0)),\rho(u(0)))$ satisfies
\begin{equation}\label{eq-cor-exi-2}
\mathfrak{m}_{\delta,\kappa}(\bm{v}_0(0),\bm{v}_1(0))^2-\norm{P_+(v_0(0)-v_1(0))}_{E_{\kappa}}^2=0<l_0^2\norm{P_+(v_0(0)-v_1(0))}_{E_\kappa}^2.
\end{equation}
Therefore, the conclusion follows Lemma \ref{lem-exit-pro} and the inequalities \eqref{eq-cor-exi-1} and \eqref{eq-cor-exi-2}.

\qed

%We define $\mathfrak{g}(w,c)$ by
%\[\mathfrak{g}(w,c)=w+G_+^{\delta}(w,c)+Q_c\]
%for $w \in H^1(\RTL), c >0$.
In the following corollary, we show the equality $\mathcal{M}_{cs}^{\delta}(c^*,\varepsilon)=\tilde{\mathcal{M}}_{cs}^{\delta}(c^*,\varepsilon)$.
\begin{corollary}\label{cor-eq-mcs}
Let $\delta,\kappa,l>0$. Assume the conditions $\eqref{ass-small}$ and $\eqref{ass-exit}$.
There exists $\varepsilon_1=\varepsilon_1(c^*,\delta)>0$ such that for $0<\varepsilon<\varepsilon_1$
\[\mathcal{M}_{cs}^{\delta}(c^*,\varepsilon)=\tilde{\mathcal{M}}_{cs}^{\delta}(c^*,\varepsilon)\]
and
\begin{align*}
&\{\mathfrak{g}(w,c^*); w \in P_{\leq 0}H^1(\RTL),  \norm{(P_1+P_2)w}_{H^1}<\varepsilon^{1/2},   \inf_{q\in \R}\norm{\mathfrak{g}(w,c^*)-\tau_qQ_{c^*}}_{H^1}<\varepsilon\}\\
=&\{\tau_{\rho}\mathfrak{g}(w,c); w \in (P_-+P_a+P_{\gamma})H^1(\RTL), |c-c^*|\leq c^*/2, \rho \in \R,\\
\ &\norm{(P_1+P_2)(\tau_{\rho}\mathfrak{g}(w,c)-Q_{c^*})}_{H^1}<\varepsilon^{1/2}, \inf_{q \in \R}\norm{\mathfrak{g}(w,c)-\tau_{q}Q_{c^*}}_{H^1}<\varepsilon\}.
\end{align*}
Moreover,
\begin{align}\label{cor-eq-mcs-1}
\mathcal{M}_{cs}^{\delta}&(c^*,\varepsilon)=\bigcup_{\rho \in \R}\{\tau_{\rho}\mathfrak{g}(w,c^*);w\in P_{\leq 0}H^1(\RTL),\notag\\
&\norm{(P_1+P_2)w}_{H^1}<\varepsilon^{1/2},\inf_{q\in \R}\norm{\mathfrak{g}(w,c^*)-\tau_qQ_{c^*}}_{H^1}<\varepsilon\}.
\end{align}

\end{corollary}
\proof
By the definitions of $\mathcal{M}_{cs}^{\delta}(c^*,\varepsilon)$ and $\tilde{\mathcal{M}}_{cs}^{\delta}(c^*,\varepsilon)$, we have
\[\mathcal{M}_{cs}^{\delta}(c^*,\varepsilon) \subset \tilde{\mathcal{M}}_{cs}^{\delta}(c^*,\varepsilon).\]
Let $0<\varepsilon<\min\{\tilde{\varepsilon}(c^*,\varepsilon^*),\varepsilon^*\}$, where $\tilde{\varepsilon}$ is defined in Theorem \ref{thm-stability} and $\varepsilon^*$ is defined in Corollary \ref{cor-exit}.
For any solutions $u(t)$ with initial data $u(0) \in \tilde{\mathcal{M}}_{cs}^{\delta}(c^*,\varepsilon)$, we have
\[\sup_{t \geq 0}\inf_{q \in \R}\norm{u(t)-\tau_qQ_{c^*}}_{H^1}<\varepsilon^*.\]
By Corollary \ref{cor-exit},  solutions $u(t)$ of \eqref{ZKeq} with initial data $u(0) \in \mathcal{N}_{\varepsilon, c^*}\setminus \mathcal{M}_{cs}^{\delta}(c^*,\varepsilon) \subset \mathcal{N}_{\varepsilon^*,c^*}\setminus \mathcal{M}_{cs}^{\delta}(c^*,\varepsilon^*)$ satisfy
\[\sup_{t \geq 0}\inf_{q \in \R}\norm{u(t)-\tau_qQ_{c^*}}_{H^1}\geq \varepsilon^*.\]
Therefore, we obtain $\tilde{\mathcal{M}}_{cs}^{\delta}(c^*,\varepsilon) \subset \mathcal{M}_{cs}^{\delta}(c^*,\varepsilon)$.

Since $\tilde{\mathcal{M}}_{cs}^{\delta}(c^*,\varepsilon) \subset \mathcal{M}_{cs}^{\delta}(c^*,\varepsilon)$, we have
\begin{align*}
&\{\mathfrak{g}(w,c^*); w \in P_{\leq 0}H^1(\RTL), \norm{(P_1+P_2)w}_{H^1}<\varepsilon^{1/2},  \inf_{q\in \R}\norm{\mathfrak{g}(w,c^*)-\tau_qQ_{c^*}}_{H^1}<\varepsilon\}\\
\subset &\{\tau_{\rho}\mathfrak{g}(w,c); w \in (P_-+P_a+P_{\gamma})H^1(\RTL), |c-c^*|\leq c^*/2, \rho \in \R,\\
\ &\norm{(P_1+P_2)(\tau_{\rho}\mathfrak{g}(w,c)-Q_{c^*})}_{H^1}<\varepsilon^{1/2}, \inf_{q \in \R}\norm{\mathfrak{g}(w,c)-\tau_{q}Q_{c^*}}_{H^1}<\varepsilon\}.
\end{align*}
Let $ w \in (P_-+P_a+P_{\gamma})H^1(\RTL)$, $c>0$ and $\rho \in \R$ satisfying $ |c-c^*|\leq c^*/2,$
\[\norm{(P_1+P_2)(\tau_{\rho}\mathfrak{g}(w,c)-Q_{c^*})}_{H^1}<\varepsilon^{1/2} \mbox{ and } \inf_{q \in \R}\norm{\mathfrak{g}(w,c)-\tau_{q}Q_{c^*}}_{H^1}<\varepsilon.
\]
We define the solution $u_0(t)$ to the equation \eqref{ZKeq} with the initial data $\tau_{\rho}\mathfrak{g}(w,c)$ and the solution $u_1(t)$ to the equation \eqref{ZKeq} with the initial data $\mathfrak{g}(P_{\leq 0 }w_0,c^*)$, where $w_0=\tau_{\rho}\mathfrak{g}(w,c)-Q_{c^*}$.
By Theorem \ref{thm-stability}, there exists $l_0>0$ such that $l_0$ and $\delta$ satisfy \eqref{ass-exit} and
\[\sup_{t \geq 0} \inf_{q \in \R}\norm{u_j(t)-\tau_qQ_{c^*}}_{H^1}<\varepsilon_*(c^*,\delta,l_0)\]
for $j=0,1$ and sufficiently small $\varepsilon>0$, where $\varepsilon_*$ is defined in Lemma \ref{lem-exit-pro}.
Since $u_1(0)$ satisfy \eqref{eq-exi-4} in Lemma \ref{lem-exit-pro} for $0<\varepsilon<\varepsilon_*$, we have
\[l_0\norm{P_+(u_0(0)-u_1(0))}_{E_{\kappa}}\lesssim \norm{P_{\leq 0}(u_0(0)-u_1(0))}_{H^1}=0\]
and $u_0(0)=u_1(0)$.
Therefore, we obtain 
\begin{align*}
&\{\mathfrak{g}(w,c^*); w \in P_{\leq 0}H^1(\RTL),  \norm{(P_1+P_2)w}_{H^1}<\varepsilon^{1/2},  \inf_{q\in \R}\norm{\mathfrak{g}(w,c^*)-\tau_qQ_{c^*}}_{H^1}<\varepsilon\}\\
\supset &\{\tau_{\rho}\mathfrak{g}(w,c); w \in (P_-+P_a+P_{\gamma})H^1(\RTL), |c-c^*|\leq c^*/2, \rho \in \R,\\
\ &\norm{(P_1+P_2)(\tau_{\rho}\mathfrak{g}(w,c)-Q_{c^*})}_{H^1}<\varepsilon^{1/2}, \inf_{q \in \R}\norm{\mathfrak{g}(w,c)-\tau_{q}Q_{c^*}}_{H^1}<\varepsilon\}.
\end{align*}
Since
\[
\inf_{\rho \in \R}\norm{(P_1+P_2)(\tau_{\rho}\mathfrak{g}(w,c^*)-Q_{c^*})}_{H^1} \leq \inf_{q\in\R}\norm{\mathfrak{g}(w,c^*)-\tau_qQ_{c^*}}_{H^1}
\]
for $w \in P_{\leq 0}H^1(\RTL)$, we have the equality \eqref{cor-eq-mcs-1} for sufficiently small $\varepsilon>0$.

\qed

\section{The $C^1$ regularity of center stable manifolds}

In this section, we prove that $G_+^{\delta}$ is a $C^1$ function on $P_{\leq 0}H^1(\RTL)\times (0,\infty)$ by applying the argument in the section 2.3 in \cite{K N S}.
Let $\varepsilon>0$ and $\psi_0,\psi_1 \in P_{\leq 0}H^1(\RTL)$ with $\norm{\psi_0}_{H^1}<\varepsilon$.
We define a solution $(v_0,c_0,\rho)$ to the system \eqref{Leq-1}--\eqref{Leq-2} such that $v_0(0)=\psi_0+G_+^{\delta}(\psi_0,c^*)$ and $c_0(0)=c^*$.
Let $v_h$ be a solution to the equation
\begin{equation}\label{eq-vh}
v_t=\partial_x\mathbb{L}_{c^*}v+(\dot{\rho}_0-c^*)\partial_xv+2\partial_x((Q_{c^*}-Q_{c_0})v)+(\dot{\rho}_0-c_0)\partial_xQ_{c_0}-\dot{c}_0\partial_cQ_{c_0}-\partial_x(v^2)
\end{equation}
with initial data $v_h(0)=\psi_0+h\psi_1+G_+^{\delta}(\psi_0+h\psi_1,c^*)$.
If $\varepsilon>0$ is sufficiently small, then $\tau_{\rho_0}(v_0+Q_{c_0})$ and $\tau_{\rho_0}(v_h+Q_{c_0})$ are solutions to the equation \eqref{ZKeq}.
By the Lipschitz continuity of $G_+^{\delta}$, for any sequence $\{h_n\}_n$ in $\R$ with $h_n \to 0$ as $n \to \infty$ there exist a subsequence $\{h_n'\}_n \subset \{h_n\}_n$ and $\psi_* \in P_+H^1(\RTL)$ such that
\[
\frac{G_+^{\delta}(\psi_0+h_n'\psi_1,c^*)-G_+^{\delta}(\psi_0,c^*)}{h_n'} \to \psi_* \mbox{ as } n \to \infty.
\]
Let $w_h=\tau_{\rho_0}v_h$ for $h\geq 0$.
Then, $w_h$ is the solution to the equation
\[
w_t=-\partial_x\Delta w -2\partial_x(\tau_{\rho_0}Q_{c_0}w)+(\dot{\rho}_0-c_0)\tau_{\rho_0}\partial_xQ_{c_0}-\dot{c}_0\tau_{\rho_0}\partial_cQ_{c_0}-\partial_x(w^2)
\]
with $w_h(0)=v_h(0).$
By the well-posedness result of the equation \eqref{ZKeq} in \cite{M P}, we have there exists $b_0>\frac{1}{2}$ such that for $T>0$ and $\frac{1}{2}  < b<b_0$ there exists $C=C(T,b)>0$ satisfying
\begin{equation}\label{c1-eq-1}
\norm{w_0}_{X_T^{1,b}}\leq C\norm{v_0(0)}_{H^1}.
\end{equation}
We define $\xi$ as the solution to the equation
\begin{equation}\label{c1-eq-xi}
\xi_t=-\partial_x\Delta \xi -2\partial_x(\tau_{\rho_0}Q_{c_0}\xi)-2\partial_x(w_0\xi)
\end{equation}
with the initial data $\xi(0)=\psi_1+\psi_*$.
By the smoothness of the flow map of equation \eqref{ZKeq} given by \cite{M P}, we have that for $T>0$ 
\begin{equation}\label{c1-eq-2}
\norm{\frac{w_{h_n'}-w_0}{h_n'}-\xi}_{L^{\infty}((-T,T),H^1)} \to 0 \mbox{ as } n \to \infty.
\end{equation}
Let $\eta=\tau_{-\rho_0}\xi$.
Then, $\eta$ satisfies the equation
\begin{equation}\label{c1-eq-eta}
\eta_t=\partial_x\mathbb{L}_{c^*}\eta-2\partial_x((Q_{c_0}-Q_{c^*})\eta)+(\dot{\rho}_0-c^*)\partial_x\eta-2\partial_x(v_0\eta).
\end{equation}

The following lemma shows the well-posedness of \eqref{c1-eq-eta}. 
\begin{lemma}\label{lem-linear-exit-0}
Let $\kappa>0$, $b_*>1/2$ and $(v_0,c_0,\rho_0)$ be a solution to the system $\eqref{Leq-1}$--$\eqref{Leq-2}$ satisfying $\sup_{t \geq 0} (\norm{v_0(t)}_{H^1}+|c_0(t)-c^*|)\leq \kappa $ and $\tau_{\rho_0(t)}(v_0(t)+Q_{c_0(t)})$ is a solution to the equation $\eqref{ZKeq}$ on $t \in [0,\infty)$ with 
\[\tau_{\rho_0(t)}(v_0(t)+Q_{c_0(t)}) \in X_{0,T}^{1,b_*}.\]
Then, the Cauchy problem of the equation $\eqref{c1-eq-eta}$ is well-posed in $H^1(\RTL)$.
Precisely, there exists $b>\frac{1}{2}$ such that for any $\eta_0 \in H^1(\RTL)$ there exists a unique solution $\xi$ to the equation $\eqref{c1-eq-xi}$ satisfying that $\xi(0)=\tau_{-\rho_0(0)}\eta_0$,
\[ \xi \in X_{0,T}^{1,b} \mbox{ for } T>0\]
and $\tau_{-\rho_0}\xi$ is a solution to the equation $\eqref{c1-eq-eta}$ with initial data $\eta_0$.
Then, the solution map $\eta_0 \to \xi \in X_{0,T}^{1,b}$ is continuous for $T>0$.
Moreover, for any solutions $\eta$ to the equation \eqref{c1-eq-eta} and $s\geq 0$ there exists $\xi_s \in X^{1,b}$ such that for $t \in (\min\{s-1,0\},s+1)$ we have $\eta(t)=\tau_{\rho_0(s)-\rho_0(t)}\xi_s(t)$ and
\begin{equation}\label{eq-c1r-4}
\norm{\eta}_{L^{\infty}((\min\{s-1,0\},s+1),H^1)} \lesssim \norm{\xi_s}_{X^{1,b}}\lesssim \norm{\eta(s)}_{H^1}.
\end{equation}
\end{lemma}

Applying the estimate \eqref{eq-w-gwp-2}--\eqref{eq-w-gwp-4} to Duhamel's formula for the equation \eqref{c1-eq-xi}, we obtain the well-posedness of the equation \eqref{c1-eq-eta}. 
We omit the detail of the proof of Lemma \ref{lem-linear-exit-0}.

In the following lemma, we prove the criterion of the growth estimate of solutions to \eqref{c1-eq-eta}.
\begin{lemma}\label{lem-linear-exit}
Let $K_0>0$.
There exist $\iota_0=\iota_0(K_0),K_1=K_1(K_0)>0$ such that for $0<\iota <\iota_0$ and a solution $(v_0,c_0,\rho_0)$ to the system $\eqref{Leq-1}$--$\eqref{Leq-2}$ satisfying that $\sup_{t \geq 0} (\norm{v_0(t)}_{H^1}+|c_0(t)-c^*|)\leq \iota $ and $\tau_{\rho_0(t)}(v_0(t)+Q_{c_0(t)})$ is a solution to the equation $\eqref{ZKeq}$ on $t \in [0,\infty)$ with 
\[\tau_{\rho_0(t)}(v_0(t)+Q_{c_0(t)}) \in X_{0,T}^{1,b_*},\]
the following holds.
If a solution $\eta$ to the equation $\eqref{c1-eq-eta}$ with initial data $\eta(0) \in H^1(\RTL)$ satisfies  
\begin{equation}\label{eq-c1r-1}
K_0\iota^{\frac{1}{3}}\norm{P_{\leq 0}\eta(t_0)}_{E_{\iota^{1/3}}}<\norm{P_+\eta(t_0)}_{E_{\iota^{1/3}}}
\end{equation}
at some $t_0\geq 0$, then for $t \geq t_0+1/2$
\begin{equation}\label{eq-c1r-2}
3\norm{P_+\eta(t)}_{E_{\iota^{1/3}}}> e^{\frac{k_*(t-t_0)}{2}}(\norm{P_+\eta(t_0)}_{E_{\iota^{1/3}}}+K_0\iota^{1/3}\norm{P_{\leq 0}\eta(t)}_{E_{\iota^{1/3}}}).
\end{equation}
On the other hand, if $\eqref{eq-c1r-1}$ fails for $t_0\geq 0$, then for $t \geq 0$
\begin{equation}\label{eq-c1r-3}
\norm{P_+\eta(t)}_{E_{\iota^{1/3}}} \leq K_0\iota^{1/3}\norm{P_{\leq 0} \eta(t)}_{E_{\iota^{1/3}}} \lesssim e^{K_1\iota^{1/6} t}\iota^{1/3}\norm{P_{\leq 0}\eta(0)}_{E_{\iota^{1/3}}}.
\end{equation}
\end{lemma}
\proof
%Applying the estimate \eqref{eq-w-gwp-2}--\eqref{eq-w-gwp-4} to Duhamel's formula for the equation \eqref{c1-eq-xi}, we obtain the well-posedness of the equation \eqref{c1-eq-eta}. 
%Moreover, for any solutions $\eta$ to the equation \eqref{c1-eq-eta} and $s\geq 0$ there exists $\xi_s \in X^{1,b}$ such that for $t \in (\min\{s-1,0\},s+1)$ we have $\eta(t)=\tau_{\rho_0(s)-\rho_0(t)}\xi_s(t)$ and
%\begin{align}\label{eq-c1r-4}
%\norm{\eta}_{L^{\infty}((\min\{s-1,0\},s+1),H^1)} \lesssim \norm{\xi_s}_{X^{1,b}}\lesssim \norm{\eta(s)}_{H^1}
%\end{align}
By the inequality \eqref{eq-w-gwp-4}, \eqref{c1-eq-1} and \eqref{eq-c1r-4}, for $t_1, t_2 \geq 0$ with $|t_1-t_2|<1$ we have
\begin{align}\label{eq-c1r-5}
|\norm{P_{\gamma}\eta(t_2)}_{E_{\iota^{1/3}}}^2-\norm{P_{\gamma}\eta(t_1)}_{E_{\iota^{1/3}}}^2| &\lesssim \iota \norm{\eta(t_1)}_{H^1}^2 \notag \\
& \lesssim \iota^{1/3}\norm{P_{\leq 0}\eta(t_1)}_{E_{\iota^{1/3}}}^2+\iota \norm{P_+\eta(t_1)}_{E_{\iota^{1/3}}}^2.
\end{align}
Since
\[
\norm{P_-\partial_t\eta(t)-\partial_x\mathbb{L}_{c^*}P_-\eta(t)}_{E_{\iota^{1/3}}}+|(Q_{c^*},\partial_t\eta(t))_{L^2}|\lesssim \iota \norm{\eta(t)}_{H^1},
\]
we have
\begin{equation}\label{eq-c1r-6}
|\norm{P_2\eta(t_2)}_{E_{\iota^{1/3}}}-\norm{P_2\eta(t_1)}_{E_{\iota^{1/3}}}| \lesssim \iota^{2/3}\norm{\eta(t_1)}_{E_{\iota^{1/3}}}
\end{equation}
and
\begin{equation}\label{eq-c1r-7}
\norm{P_-\eta(t_2)}_{E_{\iota^{1/3}}}-e^{-k_*(t_2-t_1)}\norm{P_-\eta(t_1)}_{E_{\iota^{1/3}}} \lesssim \iota^{2/3}\norm{\eta(t_1)}_{E_{\iota^{1/3}}}
\end{equation}
for $t_1,t_2 \geq 0$ with $|t_1-t_2|<1$.
By the inequality
\begin{align*}
&|(\partial_xQ_{c^*},\partial_t\eta(t))_{L^2}|+\Bigl|\Bigl(Q_{c^*}^{3/2}\cos \frac{n_0 y}{L},\partial_t \eta(t)\Bigr)_{L^2}\Bigr|+\Bigl|\Bigl(Q_{c^*}^{3/2}\sin \frac{n_0y}{L},\partial_t \eta (t)\Bigr)_{L^2}\Bigr|\\
& \lesssim \norm{P_{\gamma}\eta(t)}_{H^1}+\norm{P_2\eta(t)}_{H^1}+\iota \norm{\eta(t)}_{H^1}
\end{align*}
and \eqref{eq-c1r-4}--\eqref{eq-c1r-6} we obtain for $t_1,t_2 \geq 0$ with $|t_1-t_2|<1$
\begin{align}\label{eq-c1r-8}
&|\norm{P_1\eta(t_2)}_{E_{\iota^{1/3}}}-\norm{P_1\eta(t_1)}_{E_{\iota^{1/3}}}|+|\norm{P_a\eta(t_2)}_{E_{\iota^{1/3}}}-\norm{P_a\eta(t_1)}_{E_{\iota^{1/3}}}| \notag \\
\lesssim & \iota^{1/3}\norm{(P_{\gamma}+P_0)\eta(t_1)}_{E_{\iota^{1/3}}}+\iota^{4/3}\norm{(P_-+P_+)\eta(t_1)}_{E_{\iota^{1/3}}}
\end{align}
for sufficiently small $\iota>0$.
By the inequality \eqref{eq-c1r-5}--\eqref{eq-c1r-8}, there exists $C>0$ such that 
\begin{equation}\label{eq-c1r-9}
\norm{P_{\leq 0} \eta(t_2)}_{E_{\iota^{1/3}}}\leq (1+C\iota^{1/6})\norm{P_{\leq 0}\eta(t_1)}_{E_{\iota^{1/3}}}+C\iota^{1/2}\norm{P_+\eta(t_1)}_{E_{\iota^{1/3}}}
\end{equation}
for $t_1,t_2 \geq 0$ with $|t_1-t_2|<1$.
The inequality
\begin{equation}\label{eq-c1r-10}
\norm{P_+\partial_t\eta(t)-\partial_x\mathbb{L}_{c^*}P_+\eta(t)}_{E_{\iota^{1/3}}} \lesssim \iota \norm{\eta(t)}_{H^1}
\end{equation}
implies that there exists $C>0$ such that for $t_1,t_2\geq 0$ with $|t_1-t_2|<1$
\begin{equation}\label{eq-c1r-11}
\partial_t \norm{P_+\eta(t_2)}_{E_{\iota^{1/3}}} \geq k_* \norm{P_+\eta(t_2)}_{E_{\iota^{1/3}}}-C\iota^{2/3}\norm{\eta(t_1)}_{E_{\iota^{1/3}}}.
\end{equation}
Suppose \eqref{eq-c1r-1} for some $t_0$.
The inequality \eqref{eq-c1r-11} yields 
\begin{align}\label{eq-c1r-12}
\norm{P_+\eta(t)}_{E_{\iota^{1/3}}} \geq& e^{k_*(t-t_0)}\norm{P_+\eta(t_0)}_{E_{\iota^{1/3}}}-(e^{k_*(t-t_0)}-1)C\iota^{1/3}\norm{P_+\eta(t_0)}_{E_{\iota^{1/3}}}\notag \\
\geq & (1-C\iota^{1/3})e^{k_*(t-t_0)}\norm{P_+\eta(t_0)}_{E_{\iota^{1/3}}}
\end{align}
for $t_0\leq t<t_0+1$.
From the assumption \eqref{eq-c1r-1} and the inequalities \eqref{eq-c1r-9} and \eqref{eq-c1r-12}, we obtain
\begin{equation}\label{eq-c1r-13}
\norm{P_+\eta(t)}_{E_{\iota^{1/3}}}> (1-2C\iota^{1/3})(1+C\iota^{1/6})^{-1}e^{k_*(t-t_0)}K_0\iota^{1/3} \norm{P_{\leq 0 }\eta(t)}_{E_{\iota^{1/3}}}
\end{equation}
for $t_0 \leq t < t_0+1$ and sufficiently small $\iota>0$.
Therefore, we have
\[
\norm{P_+\eta(t)}_{E_{\iota^{1/3}}}>K_0\iota^{1/3}\norm{P_{\leq 0}\eta(t)}_{E_{\iota^{1/3}}}
\]
for $t_0+1/2 \leq t < t_0+1$ and sufficiently small $\iota>0$.
Applying this manner repeatedly, we obtain the inequality \eqref{eq-c1r-2} for $t\geq t_0+1/2$.

Suppose \eqref{eq-c1r-1} fails for $t\geq 0$.
Then, the inequality \eqref{eq-c1r-9} yields the inequality \eqref{eq-c1r-3} for all $t\geq 0$ and some $K_1>0$.
\qed

In the following lemma, we prove the uniqueness of unstable mode of the solution to \eqref{c1-eq-eta} not satisfying the growth condition \eqref{eq-c1r-1}.
\begin{lemma}\label{lem-linear-uniq}
Let $K_0>0$.
Then, there exists $\iota_1>0$ such that for $0<\iota<\iota_1$, $a_1\in \R$ and a solution $(v_0,c_0,\rho_0)$ to the system $\eqref{Leq-1}$--$\eqref{Leq-2}$ with $\sup_{t \geq 0}(\norm{v_0(t)}_{H^1}+|c_0(t)-c^*|)\leq \iota$ and for the solutions $\eta_1$ and $\eta_2$ to the equation $\eqref{c1-eq-eta}$ with $P_{\leq 0} \eta_1(0)=P_{\leq 0}\eta_2(0)$ not satisfying that $\eqref{eq-c1r-1}$ for some $t \geq 0$, we have $P_+\eta_1(0)=P_+\eta_2(0)$.
\end{lemma}

\proof
Assume there exist $0<\iota\ll \iota_0(K_0)$, a solution $(v_0,c_0,\rho_0)$ to the system \eqref{Leq-1}--\eqref{Leq-2} with $\sup_{t \geq 0}(\norm{v_0(t)}_{H^1}+|c_0(t)-c^*|)\leq \iota$ and for the solutions $\eta_1$ and $\eta_2$ to the equation \eqref{c1-eq-eta} such that $P_{\leq 0}\eta_1(0)=P_{\leq 0}\eta_2(0)$, $P_+\eta_1(0)\neq P_+\eta_2(0)$ and $\eta_1$ and $\eta_2$ do not satisfy that \eqref{eq-c1r-1} for some $t\geq 0$.
Let $\eta=\eta_1-\eta_2$.
Since 
\[\iota^{1/3}K_0\norm{P_{\leq 0}\eta(0)}_{E_{\iota^{1/3}}}< \norm{P_+\eta(0)}_{E_{\iota^{1/3}}},\]
by Lemma \ref{lem-linear-exit} we have $t\geq 1/2$
\begin{equation}\label{eq-linear-uniq-1}
3\norm{P_+\eta(t)}_{E_{\iota^{1/3}}}\geq e^{\frac{k_*t}{2}}(\norm{P_+\eta(0)}_{E_{\iota^{1/3}}}+\iota^{1/3}\norm{P_{\leq 0}\eta(t)}_{E_{\iota^{1/3}}}).
\end{equation}
On the other hand, by \eqref{eq-c1r-1} we have for $t\geq 0$
\[\norm{P_+\eta(t)}_{E_{\iota^{1/3}}}\lesssim e^{K_1\iota^{1/6}t}(\norm{P_{\leq 0}\eta_1(0)}_{E_{\iota^{1/3}}}+\norm{P_{\leq 0}\eta_2(0)}_{E_{\iota^{1/3}}})\]
and $K_1\iota^{1/6}\ll k_*$, where $K_1=K_1(K_0)$ is defined in Lemma \ref{lem-linear-exit}.
This inequality contradicts the inequality \eqref{eq-linear-uniq-1}.
Therefore, $P_+\eta_1(0)=P_+\eta_2(0)$.

\qed

We show the G\^ateaux differentiability of $G_+^{\delta}$.
Let $\delta,\kappa,l,l_0>0$ satisfying \eqref{ass-small} and \eqref{ass-exit}.
Since
\[\norm{v_h(0)}_{H^1}\lesssim \varepsilon+h\norm{\psi_1}_{H^1},\]
 Theorem \ref{thm-stability} yields
\begin{equation}\label{eq-c1rg-1}
\sup_{t\geq 0} \inf_{q \in \R}\norm{\tau_{\rho_0(t)}(v_h(t)+Q_{c_h(t)})-\tau_qQ_{c^*}}_{H^1}<\min\{ \varepsilon_*,\varepsilon_*^{1/2},\delta \}
\end{equation}
for sufficiently small $\varepsilon,h>0$, where $\varepsilon_*=\varepsilon_*(c^*,\delta,l_0,\kappa)$ is defined in Lemma \ref{lem-exit-pro}.
By $|(v_0,\partial_xQ_{c^*})_{L^2}|+|(v_0,Q_{c^*})_{L^2}|\lesssim \varepsilon_*^{1/2}$. applying Lemma \ref{lem-exit-pro}, we obtain
\[
\mathfrak{m}_{\delta,\kappa}((P_{\leq 0}v_0(t),c_0(t)),(P_{\leq 0}v_h(t),c_h(t)))>l_0\norm{P_+(v_0(t)-v_h(t))}_{E_{\kappa}}
\]
for $t \geq 0$.
Therefore, we have
\begin{equation}\label{eq-c1rg-2}
\frac{\norm{P_+(v_h(t)-v_0(t))}_{E_{\kappa}}}{\norm{P_{\leq 0}(v_h(t)-v_0(t))}_{E_{\kappa}}}\leq \frac{l_0\mathfrak{m}_{\delta,\kappa}((P_{\leq 0}v_0(t),c_0(t)),(P_{\leq 0}v_h(t),c_0(t)))}{\norm{P_{\leq 0}(v_h(t)-v_0(t))}_{E_{\kappa}}}\lesssim l_0
\end{equation}
for $t \geq 0$.
On the other hand, the convergence \eqref{c1-eq-2} yields 
\begin{equation}\label{eq-c1rg-3}
\frac{\norm{(h_n')^{-1}P_+(v_{h_n'}(t)-v_0(t))}_{E_{\kappa}}}{\norm{(h_n')^{-1}P_{\leq 0}(v_{h_n'}(t)-v_0(t))}_{E_{\kappa}}} \to \frac{\norm{P_+\eta(t)}_{E_{\kappa}}}{\norm{P_{\leq 0}\eta(t)}_{E_{\kappa}}}
\end{equation}
as $n \to \infty$ for $t \geq 0$.
By the inequality \eqref{eq-c1rg-2} and the convergence \eqref{eq-c1rg-3}, we have 
\[\frac{\norm{P_+\eta(t)}_{E_{\kappa}}}{\norm{P_{\leq 0}\eta(t)}_{E_{\kappa}}} \lesssim l_0\]
which show that the inequality \eqref{eq-c1r-2} fails for sufficiently large $t\geq 0$.
Thus,  Lemma \ref{lem-linear-exit} yields $\eta$ does not satisfy \eqref{eq-c1r-1} for $t\geq 0$.
By Lemma \ref{lem-linear-uniq}, we obtain the convergence
\begin{equation}\label{eq-c1rg-4}
\frac{G_+^{\delta}(\psi_0+h\psi_1,c^*)-G_+^{\delta}(\psi_0,c^*)}{h} \to \psi_* \mbox{ as } h \to 0
\end{equation}
which prove $G_+^{\delta}$ is G\^ateaux differentiable at $(\psi_0,c^*)$.
The linearity of the G\^ateaux derivative of $G_+^{\delta}$ follows the linearity of solutions to the equation \eqref{c1-eq-eta} with respect to the initial data.
The boundedness of the G\^ateaux derivative of $G_+^{\delta}$ follows the Lipschitz property of $G_+^{\delta}$.

Next, we prove the continuity of the G\^ateaux derivative of $G_+^{\delta}$.
Let
\[0<\varepsilon\ll \min\{\delta,\iota_0(1),K_1(1)^{-6}\}\]
and $\{\psi_n\}_{n=0}^{\infty} \subset P_{\leq 0}H^1(\RTL)$ with $\psi_n \to \psi_0$ in $H^1(\RTL)$ as $n \to \infty$ and $\sup_{n\in \Z_{\geq 0}} \norm{\psi_n}_{H^1} <\tilde{\varepsilon}(c^*,\varepsilon)$, where $\Z_{\geq 0}$ is the set of non-negative integers, $\tilde{\varepsilon}$ is defined in Theorem \ref{thm-stability} and $\iota_0$ and $K_1$ are defined in Lemma \ref{lem-linear-exit}.
By Theorem \ref{thm-stability}, we obtain 
\[\sup_{t\geq 0, n \in \Z_{\geq 0}}(\norm{v_n(t)}_{H^1}+|c_n(t)-c^*|)<\varepsilon,\]
where $(v_n,c_n,\rho_n)$ is the solution to the system \eqref{Leq-1}--\eqref{Leq-2} with initial data $(v_n(0),c_n(0),\rho_n(0))=(\psi_n,c^*,0)$.
We define $\eta_n^{\psi}$ as the solution to the equation
\begin{equation}\label{eq-c1rg-5}
\partial_t \eta=\partial_x \mathbb{L}_{c^*}\eta -2\partial_x((Q_{c_n}-Q_{c^*})\eta)+(\dot{\rho}-c^*)\partial_x\eta -2\partial_x(v_n\eta)
\end{equation}
with initial data $\psi$.
By the convergence of $\{(\tau_{-\rho_n}v_n,c_n,\rho_n)\}_n$ local in time, for $T,C>0$ we obtain the convergence 
\[\norm{\eta_n^{\psi}-\eta_0^{\psi}}_{L^{\infty}((0, T),H^1)} \to 0\]
as $n \to \infty$ uniformly on $\{\psi \in H^1(\RTL): \norm{\psi}_{H^1}<C\}$.
For $T>0$, by the boundedness of $\{\norm{\tau_{-\rho_n}v_n}_{X_T^{1,b}}+\norm{c_n-c^*}_{L^{\infty}(0,T)}\}_n$, we have the convergence
\begin{equation}\label{eq-c1rg-6}
\sup_{n\geq 0} \norm{\eta_n^{\varphi}-\eta_n^{\psi}}_{L^{\infty}((0,T), H^1)} \to 0
\end{equation}
as $\varphi \to \psi$ in $H^1$.
Let the G\^ateaux derivative of $G_+^{\delta}$ at $(\psi_n,c^*)$ be $\partial G_+^{\delta,n}$.
Applying Lemma \ref{lem-linear-exit} to $\eta_0^{\psi+\partial G_+^{\delta,0}(\psi)}$, by the boundedness of $\eta_0^{\psi+\partial G_+^{\delta,0}(\psi)}$ in time we have there exists $C>1$ such that
\begin{equation}\label{eq-c1rg-7}
3\norm{P_+\eta_0^{\psi+\partial G_+^{\delta,0}(\psi)}(t)}_{E_{\varepsilon^{1/3}}}\leq3 \varepsilon^{1/3}\norm{P_{\leq 0}\eta_0^{\psi+\partial G_+^{\delta,0}(\psi)}(t)}_{E_{\varepsilon^{1/3}}}
\leq C e^{K_1\varepsilon^{1/6}t}\varepsilon^{1/3}\norm{\psi}_{E_{\varepsilon^{1/3}}}
\end{equation}
for $t>0$, $\psi \in P_{\leq 0}H^1(\RTL)$ with $\norm{\psi}_{E_{\varepsilon^{1/3}}}\leq 1$.
On the other hand, applying Lemma \ref{lem-linear-exit} to $\eta_0^{\psi_+}$ for $\psi_+ \in P_+H^1(\RTL)\setminus \{0\}$, we have
\begin{equation}\label{eq-c1rg-8}
3\norm{P_+\eta_0^{\psi_+}(t)}_{E_{\varepsilon^{1/3}}} > e^{\frac{k_*t}{2}}\norm{\psi_+}_{E_{\varepsilon^{1/3}}}+e^{\frac{k_*t}{2}}\varepsilon^{1/3}\norm{P_{\leq 0}\eta_0^{\psi_+}(t)}_{E_{\varepsilon^{1/3}}}
\end{equation}
for $t >1/2$.
Combining the inequalities \eqref{eq-c1rg-7} and \eqref{eq-c1rg-8}, we obtain
\begin{align}\label{eq-c1rg-9}
&3\norm{P_+\eta_0^{\psi+\partial G_+^{\delta,0}(\psi)+\psi_+}(t)}_{E_{\varepsilon^{1/3}}} 
\geq  3\norm{P_+\eta_0^{\psi_+}(t)}_{E_{\varepsilon^{1/3}}}-3\norm{P_+\eta_0^{\psi+\partial G_+^{\delta,0}(\psi)}(t)}_{E_{\varepsilon^{1/3}}} \notag \\
>& e^{\frac{k_*t}{2}}\norm{\psi_+}_{E_{\varepsilon^{1/3}}}-2e^{2K_1\varepsilon^{1/6}t}C\varepsilon^{1/3}\norm{\psi}_{E_{\varepsilon^{1/3}}}+e^{K_1\varepsilon^{1/6}t}\varepsilon^{1/3}\norm{P_{\leq 0}\eta_0^{\psi+\partial G_+^{\delta,0}(\psi)}(t)}_{E_{\varepsilon^{1/3}}}\notag \\
&+e^{\frac{k_*t}{2}}\varepsilon^{1/3}\norm{P_{\leq 0}\eta_0^{\psi_+}(t)}_{E_{\varepsilon^{1/3}}}
\end{align}
for $t>1/2$, $\psi_+ \in P_+H^1(\RTL)\setminus \{0\}$ and $\psi \in P_{\leq 0}H^1(\RTL)$ with $\norm{\psi}_{E_{\varepsilon^{1/3}}}\leq 1$.
Since 
\[ \eta_n^{\varphi}(t) \to \eta_0^{\varphi}(t)\] 
 as $ n \to \infty$ in $L^{\infty}((0,T),H^1)$ uniformly on $\varphi \in H^1$ with $\norm{\varphi}_{E_{\varepsilon^{1/3}}}\leq 1$ for each $T>0$, there exists $n_T>0$ such that for $n\geq n_T$ and $1/2<t\leq T$ 
\begin{align*}
6\norm{P_+\eta_n^{\psi+\partial G_+^{\delta,0}(\psi)+\psi_+}(t)}_{E_{\varepsilon^{1/3}}}\geq & e^{\frac{k_*t}{2}}(\norm{\psi_+}_{E_{\varepsilon^{1/3}}}-2e^{-\frac{k_*t}{2}+2K_1\varepsilon^{1/6}t}C\varepsilon^{1/3}\norm{\psi}_{E_{\varepsilon^{1/3}}})\\
&+e^{K_1\varepsilon^{1/6}t}\varepsilon^{1/3}\norm{P_{\leq 0}\eta_n^{\psi+\partial G_+^{\delta,0}(\psi)+\psi_+}(t)}_{E_{\varepsilon^{1/3}}}.
\end{align*}
Thus, for $\sigma>0$ there exists $T_{\sigma}>0$ such that $\eta_n^{\psi+\partial G_+^{\delta,0}(\psi)+\psi_+}$ satisfies \eqref{eq-c1r-1} with $K_0=e^{K_1\varepsilon^{1/6}T_\sigma}$ at $T_\sigma$ for $n \geq n_{T_\sigma}$, $\psi \in P_{\leq 0}H^1(\RTL)$ and $\psi_+ \in P_+H^1(\RTL)$ with $\norm{\psi}_{E_{\varepsilon^{1/3}}}\leq 1$ and $\norm{\psi_+}_{E_{\varepsilon^{1/3}}}\geq \sigma$.
Since $\eta_n^{\psi+\partial G_+^{\delta,n}(\psi)}$ does not satisfy that \eqref{eq-c1r-1} for some $t\geq 0$, for $n \geq n_{T_\sigma}$ and $\psi \in P_{\leq 0}H^1(\RTL)$ with $\norm{\psi}_{E_{\varepsilon^{1/3}}}\leq 1$ we obtain 
\[\norm{\partial G_+^{\delta,n}(\psi)-\partial G_+^{\delta,0}(\psi)}_{E_{\varepsilon^{1/3}}}<\sigma\]
and the continuity of the G\^ateaux derivative of $G_+^{\sigma}$ at $(\psi_0,c^*)$ as the operator  from $P_{\leq 0}H^1(\RTL)$ to $P_+H^1(\RTL)$.
Therefore, $G_+^\delta$ is $C^1$ class on $P_{\leq 0}H^1(\RTL)\times \{c^*\}$ in the sense of Fr\'echet differential.
By the equation \eqref{cor-eq-mcs-1} in Corollary \ref{cor-eq-mcs}, we obtain the $C^1$ regularity of the manifold containing $\mathcal{M}_{cs}(c^*,\varepsilon)$.

% Let $G_0,G_1 \in \mathscr{G}^+_{l,\delta,\kappa}$ and $T \in [-T^*,-T^*/2]$.
% We define solutions 
% \[(v_j(t),c_j(t),\rho_j(t))=U_{\delta}(t-T)(P_{\leq 0}\psi +(\mathcal{U}_{\delta}(t)G_j)(\psi,\alpha),\alpha,q)\]
% for $j \in  \{0,1\}$, $(\psi,\alpha) \in \mathcal{H}$, $t \in \R$ and $q \in \R$.
% Then, we have 

% \qed

Funding: 
The author is supported by JSPS Research Fellowships for Young Scientists under Grant 18J00947.

% This research did not receive any specific grant from funding agencies in the public, commercial, or
% not-for-profit sectors.

\section*{Acknowledgments}
The author would like to express his great appreciation to Professor Tetsu Mizumachi for encouragements.
% The author would like to thank Professor Nikolay Tzvetkov for his helpful encouragements.
The author is supported by JSPS Research Fellowships for Young Scientists under Grant 18J00947.
% The author would like to thank Professor Tetsu Mizumachi for his helpful advices.
% The author would like to thank Professor Jean-Claude Saut for giving the reference.

% Department of Mathematics

% Kyoto University

% Kyoto 606-8502

% Japan

% E-mail address: y-youhei@math.kyoto-u.ac.jp

\begin{thebibliography}{99}


% \bibitem{ }
% 	\newblock author, 
%  	\newblock \emph{title}, 
%  	\newblock journal \textbf{vol.} (year), start page--end page.


\bibitem{A P S}
	\newblock J. C. Alexander, R. L. Pego and R. L. Sachs,
 	\newblock On the transverse instability of solitary waves in the Kadomtsev-Petviashvili equation,
 	\newblock Phys. Lett. A,  \textbf{226} (1997), 187--192. 
 
\bibitem{B J}
	\newblock P. W. Bates and C. K. R. T. Jones,
 	\newblock Invariant manifolds for semilinear partial differential equations.
 	\newblock In : Kirchgraber, U, Walther, H. O. eds. Dynam. Report. Ser. Dynam. Systems Appl, Vol. 2, Wiley, Chichester, UK, 1989,  1--38. 



\bibitem{MB}
	\newblock  M. Beceanu, 
 	\newblock A critical center-stable manifold for Schr\"odinger's equation in three dimensions,
 	\newblock Comm. Pure Appl. Math., \textbf{65} (2012), 431--507. 

\bibitem{MB0}
	\newblock M. Beceanu,
 	\newblock A center-stable manifold for the energy-critical wave equation in $\R^3$ in the symmetric setting,
 	\newblock J. Hyperbolic Differ. Equ., \textbf{11} (2014), 437--476. 
 	
\bibitem{TBB}
	\newblock T. B. Benjamin, 
 	\newblock The stability of solitary waves,
 	\newblock Proc. Roy. Soc. (London) Ser. A,  \textbf{328} (1972), 153--183. 

%\bibitem{JB}
%	\newblock  Bourgain, J.:  
% 	\newblock Fourier transform restriction phenomena for certain lattice subsets and applications to nonlinear evolution equations, I,
% 	\newblock Geom. Funct. Anal., \textbf{3}, 107--156 (1993).
 	
\bibitem{TJB}
	\newblock T. J. Bridges,
 	\newblock Universal geometric conditions for the transverse instability of solitary waves,
 	\newblock Phys. Rev. Lett.,  \textbf{84} (2000), 2614--2617.

\bibitem{C P}
	\newblock A. Comech and D. E. Pelinovsky, 
 	\newblock Purely nonlinear instability of standing waves with minimal energy, 
 	\newblock   Comm. Pure Appl. Math. \textbf{56} (2003), no. 11, 1565--1607.

 \bibitem{C M P S}
	\newblock R. C\^ote, C. Mu\~noz, D. Pilod and G. Simpson, 
 	\newblock Asymptotic Stability of High-dimensional Zakharov--Kuznetsov Solitons, 
 	\newblock  Arch. Ration. Mech. Anal. \textbf{220} (2016), no. 2, 639--710. 	

\bibitem{C R}
     \newblock M. G. Crandall and P. H. Rabinowitz,
     \newblock Bifurcation from a simple eigenvalue,
     \newblock J. Funct. Anal. \textbf{6} (1971), 1083-1102.

\bibitem{AdB}
	\newblock A. de Bouard, 
 	\newblock Stability and instability of some nonlinear dispersive solitary waves in higher dimension, 
 	\newblock  Proc. Roy. Soc. Edinburgh Sect. A \textbf{126} (1996), no. 1, 89--112.
 	
\bibitem{AVF}
	\newblock A. V. Faminskii,   
 	\newblock The Cauchy problem for the Zakharov--Kuznetsov equation,
 	\newblock Translation in Differential Equations A, \textbf{31} (1995), 1002--1012.

	
	
%\bibitem{G T V}
%	\newblock   Ginibre, J.,   Tsutsumi, Y., Velo, G.:
%	\newblock On the Cauchy problem for the Zakharov system,
%	\newblock J. Funct. Anal., \textbf{151}, 384--436  (1997).
	


\bibitem{AG1}
	\newblock A.  Gr\"unrock, 
	\newblock A remark on the modified Zakharov--Kuznetsov equation in three space dimensions,
	\newblock Math. Res. Lett., \textbf{21} (2014), 127--131.
	
	
\bibitem{G H}
	\newblock A. Gr\"unrock, S. Herr,
	\newblock The Fourier restriction norm method for the Zakharov--Kuznetsov equation,
	\newblock Discrete Contin. Dyn. Syst., \textbf{34} (2014), 2061--2068. 

%\bibitem{H K} 
%     \newblock S. Herr, S. Kinoshita,
%     \newblock Subcritical well-posedness results for the Zakharov-Kuznetsov equation in dimension three and higher,
%     \newblock arXiv:2001.09047.

	
	\bibitem{J L Z 1}
	\newblock J. Jin, Z. Lin and C. Zeng, 
	\newblock Invariant manifolds of traveling waves of the 3D Gross--Pitaevskii equation in the energy space,
	\newblock Comm. Math. Phys., \textbf{364} (2018), no. 3, 981--1039.

	\bibitem{J L Z 2}
	\newblock  J. Jin, Z. Lin and C. Zeng,   
	\newblock Dynamics near the solitary waves of the supercritical gKDV equations,
      \newblock J. Differential Equations 267 (2019), no. 12, 7213--7262.

	\bibitem{MAJ}
	\newblock M. A. Johnson,  
	\newblock The transverse instability of periodic waves in Zakharov--Kuznetsov type equations,
	\newblock Stud. Appl. Math., \textbf{124} (2010), 323--345.

\bibitem{K P} 
     \newblock B. B. Kadomtsev and  V. I. Petviashvili,  % first name middle initial. and then last name.  Only the first character in the paper title is capitalized.
     \newblock On the stability of solitary waves in weakly dispersive media,
     \newblock Sov. Phys. Dokl., \textbf{15} (1970), 539--541 .

\bibitem{Ki} 
     \newblock S. Kinoshita,
     \newblock Global Well-posedness for the Cauchy problem of the Zakharov--Kuznetsov equation in 2D,
     \newblock arXiv:1905.01490.

%\bibitem{Ki 2} 
%     \newblock S. Kinoshita,
%     \newblock Well-posedness for the Cauchy problem of the modified Zakharov-Kuznetsov equation,
%     \newblock arXiv:1911.13265.

\bibitem{K S}
     \newblock  J. Krieger and  W. Schlag, :   % first name middle initial. and then last name.  Only the first character in the paper title is capitalized.
     \newblock Stable manifolds for all monic supercritical focusing nonlinear Schr\"odinger equations in one dimension,
     \newblock J. Amer. Math. Soc., \textbf{19} (2006), 815--920. 

\bibitem{K N S 0}
     \newblock J. Krieger, K. Nakanishi and W. Schlag, % first name middle initial. and then last name.  Only the first character in the paper title is capitalized.
     \newblock Threshold phenomenon for the quintic wave equation in three dimensions,
     \newblock Comm. Math. Phys., \textbf{327} (2014), 309--332.
     
\bibitem{K N S}
     \newblock J. Krieger, K. Nakanishi, W. Schlag,   % first name middle initial. and then last name.  Only the first character in the paper title is capitalized.
     \newblock Center-stable manifold of the ground state in the energy space for the critical wave equation,
     \newblock  Math. Ann., \textbf{361} (2015),  1--50.

\bibitem{K K P}
     \newblock E. Kirr, P. G. Kevrekidis and D. E. Pelinovsky, % first name middle initial. and then last name.  Only the first character in the paper title is capitalized.
     \newblock Symmetry-breaking bifurcation in the nonlinear Schr\"odinger equation with symmetric potentials,
     \newblock Comm. Math. Phys. \textbf{308} (2011), no. 3 795-844.

\bibitem{L L S}
     \newblock D. Lannes, F. Linares, J.-C. Saut,   % first name middle initial. and then last name.  Only the first character in the paper title is capitalized.
     \newblock The Cauchy problem for the Euler-Poisson system and derivation of the Zakharov--Kuznetsov equation,
     \newblock In: Cicognani, M., Colombini, F., Del Santo, D., eds.  Studies in phase space analysis with applications to PDEs,  Progr. Nonlinear Differential Equations Appl. Vol. 84,  New York:Birkhauser/Springer,  2013, 181--213. 

\bibitem{L P R T}
     \newblock F. Linares, M. Panthee, T. Robert and N. Tzvetkov,
     \newblock On the periodic Zakharov-Kuznetsov equation. 
     \newblock Discrete Contin. Dyn. Syst. 39 (2019), no. 6, 3521--3533.

\bibitem{L P 1}
     \newblock  F. Linares and A. Pastor,  % first name middle initial. and then last name.  Only the first character in the paper title is capitalized.
     \newblock Well-posedness for the two-dimensional modified Zakharov--Kuznetsov equation,
     \newblock SIAM J. Math. Anal., \textbf{41} (2009), 1323--1339.

\bibitem{L P 2}
     \newblock F.  Linares and A. Pastor,  % first name middle initial. and then last name.  Only the first character in the paper title is capitalized.
     \newblock Local and global well-posedness for the 2D generalized Zakharov--Kuznetsov equation,
     \newblock  J. Funct. Anal., \textbf{206} (2011), 1060--1085.

\bibitem{L P S}
     \newblock F. Linares, A. Pastor and J.-C. Saut, % first name middle initial. and then last name.  Only the first character in the paper title is capitalized.
     \newblock Well-posedness for the ZK equation in a cylinder and on the background of a KdV soliton,
     \newblock Comm. Partial Differential Equations, \textbf{35} (2010), 1674--1689. 
     
\bibitem{L S}
     \newblock  F. Linares and J.-C. Saut, % first name middle initial. and then last name.  Only the first character in the paper title is capitalized.
     \newblock The Cauchy problem for the 3D Zakharov--Kuznetsov equation,
     \newblock Discrete Contin. Dyn. Syst., \textbf{24} (2009), 547--565. 
     
     
\bibitem{MM}
     \newblock M. Maeda, % first name middle initial. and then last name.  Only the first character in the paper title is capitalized.
     \newblock Stability of bound states of Hamiltonian PDEs in the degenerate cases,
     \newblock J. Funct. Anal.  \textbf{263 } (2012), no. 2, 511--528.

\bibitem{M M 1}
	\newblock Y. Martel and F. Merle, 
 	\newblock Asymptotic stability of solitons for subcritical generalized KdV equations,
 	\newblock Arch. Ration. Mech. Anal., \textbf{157} (2001), 219--254.

\bibitem{M M 2}
	\newblock Y. Martel and F. Merle,
 	\newblock Asymptotic stability of solitons of the subcritical gKdV equations revisited,
 	\newblock Nonlinearity, \textbf{18} (2005), 55--80.
 	
 \bibitem{M M 3}
	\newblock Y. Martel and F. Merle, 
 	\newblock Asymptotic stability of solitons of the gKdV equations with general nonlinearity,
 	\newblock Math. Ann., \textbf{341} (2008), 391--427.

 \bibitem{M M N R}
	\newblock Y. Martel, F. Merle, K. Nakanishi and P. Rapha\"el,  
 	\newblock Codimension one threshold manifold for the critical gKdV equation,
 	\newblock Comm. Math. Phys., \textbf{342} (2016), 1075--1106.
 	
\bibitem{TM1}
	\newblock T. Mizumachi, 
 	\newblock Large time asymptotics of solutions around solitary waves to the generalized Korteweg--de Vries equations,
 	\newblock SIAM J. Math. Anal., \textbf{32} (2001), 1050--1080.
     
\bibitem{TM2}
	\newblock T. Mizumachi, 
 	\newblock Stability of line solitons for the KP-II equation in $\R^2$,
 	\newblock Mem. Amer. Math. Soc., \textbf{238} (2015), vii+95 pp. . 

\bibitem{TM3}
	\newblock T. Mizumachi,
 	\newblock Stability of line solitons for the KP-II equation in $\R^2$, II,
 	\newblock Proc. Roy. Soc. Edinburgh Sect. A, \textbf{148} (2018), 149--198. 
 	
 \bibitem{M T}
	\newblock T. Mizumachi and N. Tzvetkov,
 	\newblock Stability of the line soliton of the KP-II equation under periodic transverse perturbations,
 	\newblock Math. Ann., \textbf{352} (2012), 659--690. 
 	
\bibitem{M P}
	\newblock  L. Molinet and D. Pilod, 
	\newblock  Bilinear Strichartz estimates for the Zakharov--Kuznetsov equation and applications,
	\newblock  Ann. Inst. H. Poincar\'e Anal. Non Lineaire, \textbf{32} (2015), 347--371.

\bibitem{M S T}
    \newblock L. Molinet, J.-C. Saut and N. Tzvetkov,   
    \newblock Global well-posedness for the KP-II equation on the background of a non-localized solution,
    \newblock Ann. Inst. H. Poincar\'e Anal. Non Lineaire, \textbf{28} (2011), 653--676.

	
\bibitem{N S}
	\newblock  K. Nakanishi and W. Schlag,
	\newblock  Global dynamics above the ground state energy for the cubic NLS equation in 3D,
	\newblock  Calc. Var. Partial Differential Equations, \textbf{44} (2012), 1--45.
	
\bibitem{N S 2}
	\newblock  K. Nakanishi and W. Schlag,
	\newblock  Invariant manifolds around soliton manifolds for the nonlinear Klein--Gordon equation,
	\newblock  SIAM J. Math. Anal., \textbf{44} (2012), 1175--1210. 

%  	
\bibitem{P W 2} %(MR1177566) [10.1098/rsta.1992.0055]
\newblock  R. Pego and M. I. Weinstein,  
\newblock Asymptotic stability of solitary waves,
\newblock Comm. Math. Phys., \textbf{164} (1994), 305--349. 

\bibitem{DP} 
\newblock D. Pelinovsky, 
\newblock Normal form for transverse instability of the line soliton with a nearly critical speed of propagation,
\newblock Math. Model. Nat. Phenom., \textbf{13} (2018), 1--20.

\bibitem{R V}
	\newblock F. Ribaud and S. Vento, 
	\newblock Well-posedness results for the three-dimensional Zakharov--Kuznetsov equation,
	\newblock SIAM J. Math. Anal., \textbf{44} (2012), 2289--2304.
	
\bibitem{R T 0}
	\newblock F. Rousset and N. Tzvetkov, 
	\newblock Transverse nonlinear instability of solitary waves for some Hamiltonian PDE's,
	\newblock J. Math. Pures. Appl., \textbf{90} (2008), 550--590.

\bibitem{R T 1}
	\newblock F. Rousset and N. Tzvetkov, 
	\newblock Transverse nonlinear instability for two-dimensional dispersive models,
	\newblock Ann. I. Poincar\'e-AN, \textbf{26} (2009), 477--496.


\bibitem{R T 3} %(MR2928221) [10.1007/s00220-012-1495-y]
\newblock  F. Rousset and N. Tzvetkov, 
\newblock Stability and instability of the KdV solitary wave under the KP-I flow,
\newblock Comm. Math. Phys., \textbf{313} (2012), 155--173.

\bibitem{WS} %(MR1752509) [10.1090/conm/255/03982]
\newblock W. Schlag,
\newblock Stable manifolds for an orbitally unstable nonlinear Schr\"odinger equation,
\newblock Ann. of Math. (2), \textbf{169} (2009), 139--227. 



\bibitem{V A}
     \newblock J. Villarroel and M. J. Ablowitz,   % first name middle initial. and then last name.  Only the first character in the paper title is capitalized.
     \newblock On the initial value problem for the KPII equation with data that do not decay along a line,
     \newblock Nonlinearity, \textbf{17} (2004), 1843--1866. 

\bibitem{YY2}
     \newblock Y. Yamazaki, % first name middle initial. and then last name.  Only the first character in the paper title is capitalized.
     \newblock \emph{Stability of line standing waves near the bifurcation point for nonlinear Schr\"odinger equations},
     \newblock  Kodai Math. J. \textbf{38} (2015), no. 1, 65--96.
     
%     \bibitem{YY3}
%     \newblock  Yamazaki, Y.: % first name middle initial. and then last name.  Only the first character in the paper title is capitalized.
%     \newblock Transverse instability for nonlinear Schr\"odinger equation with a linear potential,
 %    \newblock Adv. in Differential Equations, \textbf{21}, 429--462 (2016).

     \bibitem{YY4}
     \newblock  Y. Yamazaki,  % first name middle initial. and then last name.  Only the first character in the paper title is capitalized.
     \newblock Stability for line solitary waves of Zakharov--Kuznetsov equation,
     \newblock J. Differential Equations, \textbf{262} (2017), 4336--4389.

     \bibitem{YY5}
     \newblock  Y. Yamazaki,  % first name middle initial. and then last name.  Only the first character in the paper title is capitalized.
     \newblock Center stable manifolds around line solitary waves of the Zakharov--Kuznetsov equation,
     \newblock arXiv:1808.07315.
\bibitem{VEZ}
     \newblock V. E. Zakharov,  % first name middle initial. and then last name.  Only the first character in the paper title is capitalized.
     \newblock Instability and nonlinear oscillations of solitons,
     \newblock  JETP Lett., \textbf{22} (1975), 172--173. 
     
\bibitem{Z K}
     \newblock V. E. Zakharov and E. A.  Kuznetsov,  % first name middle initial. and then last name.  Only the first character in the paper title is capitalized.
     \newblock On three dimensional solitons,
     \newblock  Sov. Phys.-JETP,  \textbf{39} (1974), 285--286.


\end{thebibliography}
\end{document}